\documentclass{siamltex}
\usepackage[utf8]{inputenc} 
\usepackage{standalone}
\usepackage{graphicx}
\usepackage{amsmath,amsfonts,amssymb}

\usepackage{tikz}
\usepackage[Publish]{externaltikz}
\usepackage{pgfplots}
\usepgfplotslibrary{groupplots}
\usepackage{subcaption}

\newcommand{\dis}{\displaystyle}

\begin{document}
\title{Kinetic and related macroscopic models for chemotaxis on networks}

\author{R. Borsche \footnotemark[1]\and A. Klar\footnotemark[1] \footnotemark[2]\and T.N.H. Pham \footnotemark[1]}
\footnotetext[1]{Technische Universit\"at Kaiserslautern, Department of Mathematics, Erwin-Schr\"odinger-Stra{\ss}e, 67663 Kaiserslautern, Germany 
  (\{borsche,klar,tpham\}@mathematik.uni-kl.de)}
\footnotetext[2]{Fraunhofer ITWM, Fraunhoferplatz 1, 67663 Kaiserslautern, Germany}

\maketitle

\begin{abstract}
In this paper we consider kinetic and associated macroscopic models for chemotaxis on a network.
Coupling conditions at the nodes of the network for the kinetic problem are presented and used to derive coupling conditions for the macroscopic approximations. 
The results of the different models are compared and relations to a Keller-Segel model on networks are discussed.
For a numerical approximation of the governing equations asymptotic preserving relaxation schemes are extended to directed graphs.
Kinetic and macroscopic equations are investigated numerically and their solutions are compared for tripod and more general networks.
\end{abstract}


\section{Introduction}

This paper focuses on  kinetic and macroscopic models for chemotaxis, describing the movement of organisms or cells following  the distribution of a chemical substance. 
In  the present context we consider  theses models  to describe chemotaxis on a network or graph.
A typical application of chemotaxis models in such geometries is given by tissue-engineering  and  the movement of fibroblasts on artificial scaffolds during  wound healing,  see \cite{refId0} for more details and further references.
 
The original equation to model chemotaxis are the  Keller-Segel equations.
These equations and in particular, the properties of their solutions  have been intensively investigated, see for example \cite{CA06,HE97,KS70,RA95}. 
There are many  adapted and expanded versions of the Keller-Segel model \cite{BO10,KU12,CH81,HI09,KS70,KS71,Chav06}.
Moreover, in recent literature improved flux-limited Keller-Segel models, taking into account the finite speed of propagation, have been developed in \cite{BBNS10}.
For a survey and an extended reference list, see for example \cite{BBNS12}.

Our starting point is the  classical kinetic chemotaxis  equation \cite{CMPS04}.
Scaling it with the so called  diffusive scaling leads to the  Keller-Segel equation, see \cite{CMPS04}. 
In general, the derivation of  Keller-Segel type models, including flux-limited diffusion models and  Fokker-Planck type  models, from underlying kinetic or microscopic models is discussed for example in  \cite{BBNS12,Chav06,CMPS04}.
In particular, using moment closure approaches one may obtain macroscopic equations
intermediate between kinetic and Keller-Segel equations, see the above mentioned references
or \cite{Filbet,Hillen}.
In particular, full moment equations with a linear closure function are sometimes called the P1-model (in the radiative transfer literature) or the Cattaneo equations.
Additionally, we introduce half-moment closures for the  kinetic equations and obtain asssoiated hydrodynamic macroscopic equations similar to \cite{dubroca}. 
Using the diffusion scaling for the half moment or  the Cattaneo equations one obtains again the Keller-Segel equations as scaling limit. 
\begin{figure}[ht!]
\centering
\includegraphics[width=0.7\textwidth]{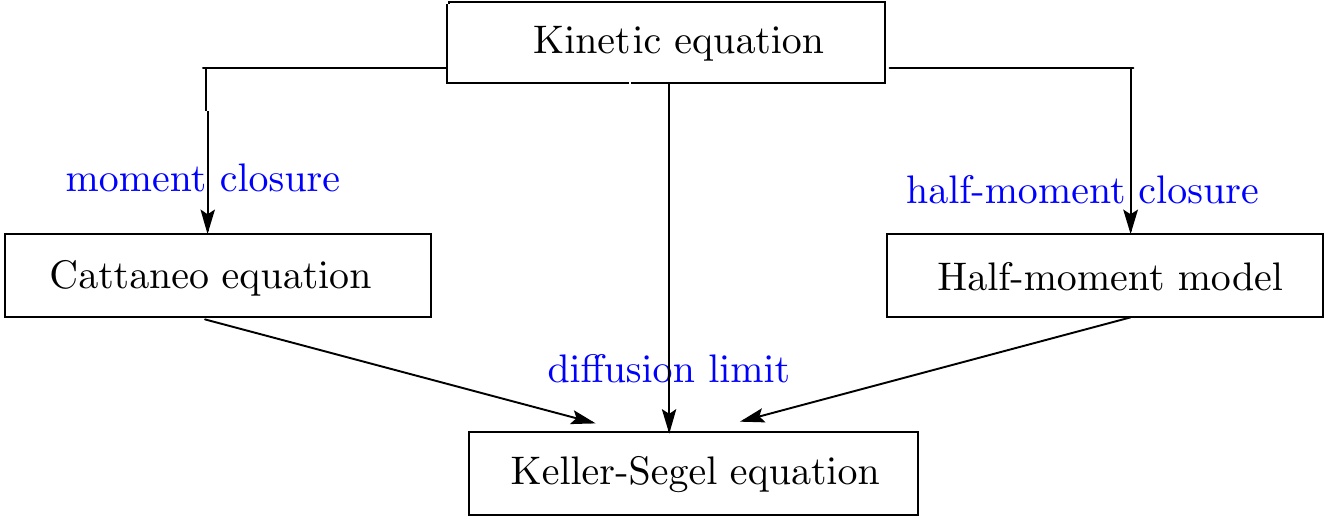}
\caption{Hierarchy of chemotaxis models}
\label{fig:Model_Limits}
\end{figure}
 
 To model chemotaxis on a network, the crucial point  is to define suitable coupling conditions.
 In previous work  coupling conditions for the scalar Keller-Segel equation \cite{MR3149314}
 have been discussed.
 Moreover, in \cite{refId0,2014arXiv1411.6109R} the Cattaneo equations on a graph are studied.
  
 Here we will develop coupling conditions for the kinetic  chemotaxis equations and derive from these equations coupling condition for the macroscopic Cattaneo  and half-moment models. 
 These  coupling conditions  guarantee on the one hand  the conservation of mass through nodes, and on the other hand, they satisfy a positivity  condition for the  density. 
 Moreover, in the diffusive limit, when the scaling parameter goes  to $0$, all coupling conditions converge to those of the Keller-Segel model, i.e. the conservation
 of mass through nodes and  a continuity condition, compare for example \cite{MR3149314}.

Numerical methods to deal with these equations have been developed in  many publications, see \cite{KU08,KU12} for the Keller-Segel equations.
For the kinetic and hydrodynamic equations, the scaling parameter describing the diffusive scaling  may differ by orders of magnitude
for different regimes. Thus, it is desirable to use asymptotic preserving schemes, a class of numerical methods that can work uniformly with respect to this parameter and in particular near to the diffusive limit, see \cite{carillo,MR1655853,refId0,Gosse,natalini,Gosse2} for the chemotaxis case.
We investigate the kinetic equations, the Cattaneo equations and the half-moment models 
on a single line using diffusive relaxation schemes as developed in \cite{MR1322811} and \cite{MR1655853}.
Combining this with a numerical scheme to solve the coupling problem we obtain  results for the full  network problem.

The paper is organized as follows. In section 2 we discuss the kinetic model
and the derivation of even-odd parity models, the Keller-Segel model and hydrodynamic model like the Cattaneo and a half moment model. Section 3 considers the coupling conditions
for the kinetic equations and the derived conditions for the macroscopic equations.
Section 4 and 5 contain details about  the numerical schemes on an interval and the numerical treatment of the coupling conditions.
Finally, the numerical results for tripod and more general networks are shown in Section 6.

\section{Chemotaxis models}
 \label{sec:models}
 
\subsection{The kinetic chemotaxis model}

Denote with   $f=f(x,t,v)$  the density of  cells at position $x\in\mathbb{R}$, 
moving with velocity $v\in V\equiv [-1,\space 1]$ at the time $t\in [0,\space T]$, 
$\rho(x,t)=\int_Vf(x,t,v)dv$ is the macroscopic density of cell,  
$m=m(x,t)$ is the density of the chemoattractant.
As in \cite{CMPS04} we have
 \begin{equation}\label{equ:g_kinetic}
  \left\{
  \begin{array}{lcl}
 \partial_t f + v\partial_x f&=& L\left[f\right] \\[0.3cm]
  \partial_t m -D\partial_{xx}m&=&\gamma_{\rho}\rho-\gamma_m m\\
  \end{array}
  \right.
 \end{equation}
with the definition
$$L[f]=L_1[f]+L_2[f].$$
Here,
$$L_i[f](x,t,v)=\int_V\left(k_i(v,v')f(x,t,v')-k_i(v',v)f(x,t,v)\right)dv'$$ for $i=1,2$.
Velocity changes are described by the probability kernel  $k_i(v,v') = k_i(x,t,v,v')$, where we neglect a possible dependence on $(x,t)$ in our notation.
New   velocities are denoted by  $v\in V $,  previous velocities b $v' \in V$.
We assume that
\begin{equation*}
 \left\{
 \begin{array}{lcl}
\int_V v'k_1(v',v)dv'&=&0\\[0.3cm]
\int_V k_1(v',v)dv'&=&\lambda\ .
 \end{array}
 \right.
 \end{equation*}
 Moreover, we assume that there exists an equilibrium distribution $F(v)$ such that the detailed balance condition:
\begin{equation}
k_1(v',v)F(v)=k_1(v,v')F(v')
\end{equation}
and
\begin{equation*}
 \left\{
 \begin{array}{lcl}
 \int_VF(v)dv&=&1 \\[0.3cm]
 \int_VvF(v)dv&=&0\ .
 \end{array}
 \right.
 \end{equation*}
For the following consideration of chemotaxis processes on networks,
we concentrate on a simple example. 
We choose $F(v)=\dis\frac{1}{2}$, 
$k_1(v,v') = \dis\frac{\lambda}{2}$ and 
$k_2(v,v')=\dis\frac{1}{2}\alpha v\overline{\partial_x m}$, with  $\overline{\partial_x m} = \frac{\partial_x m}{\sqrt{1+\vert\partial_x m\vert^2}}$.

This leads to the  the kinetic equation 
 \begin{equation}\label{equ:fi_kinetic}
 \left\{
 \begin{array}{lcl}
 \partial_t f + v\partial_x f&=&-\lambda \left( f - \frac{\rho}{2} \right)+\dis\frac{1}{2}\alpha v\dis\overline{\partial_x m}\\[0.3cm]
\partial_t m -D(\partial_{xx}m)&=&\gamma_{\rho}\rho -\gamma_m m\ . 
\end{array}
 \right.
 \end{equation}
To guarantee existence and uniqueness of the problem we assume for all $v$, as in \cite{CMPS04}
 $$\frac{\lambda}{2}+\frac{1}{2}\alpha v\overline{\partial_x m}\geq 0.$$ This 
leads to 
\begin{equation}
\label{condprev}
\lambda\geq \alpha.
\end{equation}
We introduce a parameter $\epsilon\in [0, \infty)$  and  scale the  kinetic equation \eqref{equ:fi_kinetic} in the following way 
\begin{equation}\label{equ:sca_kinetic}
  \left\{
 \begin{array}{lcl}
 \partial _t f + \dis\frac{1}{\epsilon} v\partial_x f &=& -\dis\frac{\lambda}{\epsilon^2}\left( f - \frac{\rho}{2}\right)+\frac{1}{2\epsilon}\alpha v \overline{\partial_x m}\rho\\[0.3cm]
 \partial _t m - D(\partial_{xx}) m&=& \gamma_\rho \rho-\gamma_m m \\
 \end{array}
 \right.
 \end{equation} 
Condition \eqref{condprev} turns into $\lambda\geq \epsilon \alpha$ or
 \begin{equation}
 \label{cond}
 \epsilon \le \frac{\lambda}{\alpha}.
 \end{equation}
This equation will be used as the starting point for the macroscopic models considered in the next section.

The kinetic equation \eqref{equ:g_kinetic} can  be  transformed  using the even- and odd-parities \cite{MR1655853}
\begin{align}\label{eq:evenodd}
 r(x,t,v)=\dis\frac{1}{2}\left(f(x,t,v)+f(x,t,-v)\right) ,&&
 j(x,t,v)=\dis\frac{1}{2\epsilon}\left(f(x,t,v)-f(x,t,-v)\right)
 \end{align}
 for positive velocities $(v\geq 0)$.
 This leads to the system
 \begin{align}
 \label{equ:oe evolution}
 \partial_t r + v\partial_x j&=-\dis\frac{\lambda}{\epsilon^2}\left( r-\frac{\rho}{2}\right) \\
 \partial_t j + \dis\frac{1}{\epsilon^2}v\partial_x r&=-\dis\frac{1}{\epsilon^2}\left(\lambda j -\frac{1}{2}\alpha v 
\overline{\partial_x m}
 \rho\right)\ .
 \end{align}

 When $\epsilon\to 0$, the limit of the kinetic equation \eqref{equ:sca_kinetic} is  the Keller-Segel equations with flux-limited chemoattractant.\\
 The procedure is shortly revisited.
 Start with the  scaled version of the  kinetic equation \eqref{equ:oe evolution}.
 Integrating the equation for $r$ over $v \ge 0$ we obtain
 \begin{align}
  \label{cont}
 \partial_t \rho + 2 \partial_x \int_{0}^1 v j dv =0.
  \end{align}
 Moreover we observe 
 $
 r= \frac{\rho}{2} + \mathcal{O}(\epsilon^2)
 $
 and
 $
 j = \frac{\alpha \overline{\partial_x m} }{2 \lambda} v \rho - \frac{v}{\lambda}\partial_x r 
 +\mathcal{O}(\epsilon^2) $ .
Plugging the expansions into \eqref{cont} gives the diffusion limit of the kinetic equation, the modified or flux-limited Keller-Segel equation
 \begin{equation}\label{equ:d_limit}
 \left\{
 \begin{array}{lcl}
 \partial_t\rho - \partial_x \left(\dis\frac{1}{3\lambda}\partial_x\rho - \dis\frac{\alpha}{3\lambda}\dis\left(\frac{\partial_x m}{\sqrt{1+\vert\partial_xm\vert^2}}\right)\rho\right)&=&0\\[0.3cm]
 \partial_t m -D\partial_{xx}m&=&\gamma_{\rho}\rho -\gamma_m m\\
 \end{array}.
 \right.
 \end{equation}
 
 \subsection{Macroscopic full moment model /  Cattaneo equations / P1 model}
 Similarly, we also construct the moment macroscopic models for \eqref{equ:sca_kinetic} by applying the linear moment closure.
 Therefore consider the following averaged quantities
 \begin{align*}
 \rho(x,t)=\int_{-1}^1f(x,t,v)dv\ ,&& q(x,t)=\dis\frac{1}{\epsilon}\int_{-1}^1 vf(x,t,v)dv\ ,
 \end{align*}
 where $\rho$ represents the density and $q$ the flow of the cells.
 
 Now we consider the following linear closure
 $ f(x,t,v) = a\rho(x,t) + \epsilon b vq(x,t)$ .
 Inserting this into the above averages we obtain the values $a = \frac12$ and  $b = \frac32$ such that
 \begin{equation*}
 f(x,t,v) = \dis\frac{1}{2}\rho(x,t) + \epsilon \dis\frac{3}{2} vq(x,t)\ .
 \end{equation*}
 this leads to a macroscopic model for chemotaxis as in\cite{refId0,2014arXiv1411.6109R}
 \begin{equation}\label{equ:p1}
  \left\{
  \begin{array}{lcl}
  \partial_t\rho + \partial_x q&=&0 \\
  \partial_t q +\dis\frac{1}{3\epsilon^2}\partial_x\rho & =&-\dis\frac{1}{\epsilon^2}\left( \lambda q  - \frac{\alpha}{3}\left(\frac{\partial_x m}{\sqrt{1 +\vert\partial_x m\vert^2}} \right)\rho\right)\\
  \partial_tm - D\partial_{xx}m &=& \gamma_{\rho}\rho-\gamma_mm\ .
  \end{array}
  \right.
  \end{equation}
 We note that the diffusive subcharacteristic condition, see \cite{JL}, requires $0\leq\epsilon\leq \dis\frac{\sqrt{3}\lambda}{\alpha}$, which is fulfilled by our assumption
 \eqref{cond}.
  When $\epsilon\to 0$,  the $P_1$-model (\ref{equ:p1}) has the same macroscopic diffusive limit as  the kinetic equation (\ref{equ:sca_kinetic}), i.e. the Keller-Segel equations.

\subsection{Macroscopic half-moment model}
We will construct a half-moment macroscopic model for \eqref{equ:sca_kinetic} by applying a half-moment linear closure, see  \cite{dubroca} for a nonlinear version. 
We consider the following averaged quantities
\begin{equation}
\begin{aligned}
\rho ^- &=\dis\int_{-1}^0 f(v)dv\ , \qquad
& \rho^+ &=\dis\int_{0}^1f(v)dv\ ,
\\
 q ^- &=\dis\int_{-1}^0 vf(v)dv \ ,
& q^+ &=\dis\int_{0}^1vf(v)dv \ .
\end{aligned}
\label{eq:halfmoment_average}
\end{equation}
The closure assumption is that the half moments, i.e. for positive or negative velocity, are affine linear functions in $v$ 
\begin{align}\label{eq:halfmoment_closure}
f(v) = a^+ + vb^+,\ v\geq 0&&\text{and}&&
f(v) = a^- + vb^-,\ v\leq 0
\ .
\end{align}
Inserting into \eqref{eq:halfmoment_average} leads to
\begin{align*}
\begin{aligned}
\rho^- &= a^- -\dis\frac{1}{2}b^-
\ ,\quad&
\rho^+ &=a^+ + \dis\frac{1}{2}b^+\ , \\[0.3cm]
q^- &= -\dis\frac{1}{2}a^- + \dis\frac{1}{3}b^- 
\ ,\quad &
q^+ &= \dis\frac{1}{2}a^+ + \dis\frac{1}{3}b^+\ .
\end{aligned}
\end{align*}
Furthermore we obtain
\begin{align*}
\dis\int_0^1 v^2f(v)dv = -\dis\frac{1}{6}\rho^+ + q^+ \ ,
&&
\dis\int_{-1}^0 v^2f(v)dv = -\dis\frac{1}{6}\rho^- - q^-\ .
\end{align*}
Finally, we get the half-moment system for kinetic equation \eqref{equ:sca_kinetic} as
\begin{equation}\label{equ:sc half-moment system}
\left\{
\begin{array}{lcl}
\epsilon\partial_t\rho^+ + \partial_xq^+ &=& -\dis\frac{1}{\epsilon}\lambda\left(\rho^+ -\dis\frac{\rho^+ +\rho^-}{2}\right)+\dis\frac{1}{4}\alpha \overline{\partial_x m}\left(\rho^+ +\rho^-\right)\\[0.3cm]
\epsilon\partial_tq^+ + \partial_x\left(-\dis\frac{1}{6}\rho^+ + q^+\right) &=& -\dis\frac{1}{\epsilon}\lambda\left(q^+ -\dis\frac{\rho^+ +\rho^-}{4}\right) + \dis\frac{1}{6}\alpha \overline{\partial_x m}\left(\rho^+ +\rho^-\right)\\[0.3cm]
\epsilon\partial_t\rho^- + \partial_xq^- &=& -\dis\frac{1}{\epsilon}\lambda\left(\rho^- -\dis\frac{\rho^+ +\rho^-}{2}\right)-\dis\frac{1}{4}\alpha \overline{\partial_x m}\left(\rho^+ +\rho^-\right)\\[0.3cm]
\epsilon\partial_tq^- + \partial_x\left( -\dis\frac{1}{6}\rho^- - q^- \right) &=& -\dis\frac{1}{\epsilon}\lambda\left(q^- +\dis\frac{\rho^+ +\rho^-}{4}\right) + \dis\frac{1}{6}\alpha \overline{\partial_x m}\left(\rho^+ +\rho^-\right)\\
\end{array}
\right.
\end{equation}
By introducing the variables
\begin{align}
\begin{aligned}
\rho &= \rho^+ + \rho^-\ ,
& \qquad
\hat{\rho} &= \dis\frac{\rho^+ - \rho^-}{\epsilon}\ ,
\\
q &= \dis\frac{q^+ + q^-}{\epsilon}\ , 
&
\hat{q} &= q^+ -q^-\ ,
\end{aligned}
\label{eq:halfmoment_defrhoqrhohatqhat}
\end{align}
we can rewrite the system as 
\begin{equation}\label{equ:sc half-moment system_hat}
\left\{
\begin{array}{lcl}
\partial_t\rho + \partial_xq &=& 0\\[0.3cm]
\partial_tq + \dis\frac{1}{\epsilon^2}\partial_x\left(-\dis\frac{1}{6}\rho + \hat{q}\right) &=&-\dis\frac{1}{\epsilon^2}\left(\lambda q - \dis\frac{\alpha}{3}\overline{\partial_x m}\rho\right)\\[0.3cm]
\partial_t\hat{\rho} + \dis\frac{1}{\epsilon^2}\partial_x\hat{q} &=&  -\dis\frac{1}{\epsilon^2}\left(\lambda\hat{\rho} -\dis\frac{1}{2}\alpha \overline{\partial_x m}\rho\right)\\[0.3cm]
\partial_t\hat{q} + \partial_x\left(-\dis\frac{1}{6}\hat{\rho} + q\right) &=& -\dis\frac{1}{\epsilon^2}\left(\hat{q}-\dis\frac{\rho}{2}\right)\ .
\end{array}
\right.
\end{equation}
The  half-moment model \eqref{equ:sc half-moment system} has again the Keller-Segel equations
as  macroscopic diffusive limit as $\epsilon$ goes to $0$.

\section{Coupling conditions for chemotaxis models on networks} 
\label{sec:CouplingConditions}
In this section we develop a hierarchy of coupling conditions for the models proposed in the previous section.
We start with the construction of coupling conditions for the kinetic model.
These can be used to derive coupling conditions of the remaining models.
In case of the Cattaneo model we also recall alternative coupling conditions, e.g. as proposed in \cite{2014arXiv1411.6109R}. 
In section \ref{sec:numerical_coupling} a numerical method for the coupling conditions is discussed.

In general networks are complicated structures. 
Nevertheless, for the development of coupling conditions it is sufficient to investigate a single node.
Once coupling conditions for one node are specified, arbitrary networks can be constructed.
Therefore we will discuss only the case of a single node.
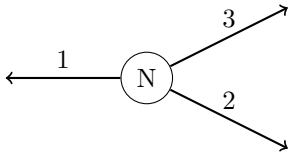
\begin{figure}[ht!]
\centering
\externaltikz{tripod_net}{

\begin{tikzpicture}
[
ls/.style={line width=2pt}]
\def\dy{1}
\def\dx{2}


\node (a) at (0,0) {};
\node[shape=circle,draw=black] (b) at (\dx,0) {N};
\node (c) at (2*\dx,\dy) {};
\node (d) at (2*\dx,-\dy) {};


\draw[<-,thick] (a)--(b)node[above,pos=0.5]{1};
\draw[->,thick] (b)--(c)node[above,pos=0.5]{3};
\draw[->,thick] (b)--(d)node[above,pos=0.5]{2};


\end{tikzpicture}
}
\caption{Sketch of a tripod network}
\label{fig:out_tripod}
\end{figure}

\subsection{Coupling condition for kinetic equations}\label{sec:cc_kinetic}
The coupled system \eqref{equ:sca_kinetic} is composed of two equations, 
a linear hyperbolic part
\begin{equation}\label{eq:kinetic_f}
\partial_t f +\dis\frac{1}{\epsilon} v\partial_x f = -\dis\frac{\lambda}{\epsilon^2}\left( f   - \dis\frac{\rho}{2} \right)+\dis\frac{1}{2 \epsilon}\alpha v
\overline{\partial_x m} \rho
\end{equation}
and a linear parabolic one for $m$.
Since these two equations are coupled via the source terms on the right hand sides, 
we can develop coupling conditions for both equations separately. 
Appropriate coupling conditions for the parabolic part can be found e.g. in \cite{MR3149314} and will not be addressed in the following.
Alternative transmission conditions are studied in \cite{2014arXiv1411.6109R}.
In contrast for \eqref{eq:kinetic_f} general coupling conditions are not known,
only special cases are considered e.g. in \cite{Herty_kinetic}.

In the following we will develop coupling conditions for a junction which connects edges governed by \eqref{eq:kinetic_f}.
For hyperbolic equations the number of coupling conditions only depends on the number of characteristics leaving the node \cite{MR2377285,MR2592880}.
Thus in the case of \eqref{eq:kinetic_f} the coupling conditions should assign a value to all $f(v)$ with $v\in[0,1]$.
In the construction of the coupling conditions we impose the following conditions:
\begin{enumerate}
 \item The coupling conditions should be linear.
 \item The coupling conditions should be independent of $v$.
 \item The total mass in the system should be conserved.
 \item The values of $f$ should remain positive all times.
 \item In the case of a $1$-to-$1$ coupling the solution of the coupled problem should coincide with the solution on one continuous edge.
 \item In the limit $\epsilon\rightarrow 0$ the conditions should converge to the coupling conditions for the Keller-Segel equation, i.e. the continuity of the densities $\rho_i$.
\end{enumerate}
Especially the first and second condition are chosen to simplify the computations.
If needed for a more detailed modeling, these can be weakened accordingly, which forces some slight modifications in the considerations below. 
In the following we will discuss the case of a junction connecting $N$ edges, 
which are all oriented away from the node.
This can be extended to junctions of an arbitrary number of in- and outgoing edges by simple local transformations.

From the first two requirements we can conclude that the coupling conditions are given in the form
\begin{align*}
 f^+ = A f^-
\end{align*}
where $f_i^+=f_i(v)$ and $f_i^-=f_i(-v)$ for $v\in[0,1]$ and $i = 1,\dots,N$.
In order to conserve the total mass in the system the matrix $A\in \mathbb{R}^{N\times N}$ has to fulfill
\begin{align*}
 \sum_{i=1}^{N}a_{i,j} = 1 \qquad \forall j = 1,\dots,N\ .
\end{align*}
The requirement 4. forces all entries of the matrix to be positive $a_{i,j}\geq 0$, $i,j=1,\dots,N$.
In the case of a $1$-to-$1$ coupling, i.e. $N=2$, we want 
\begin{align*}
 A
 =\left(\begin{array}{cc}
  0&1\\1&0
  \end{array}\right)
  \ .
\end{align*}
This we generalize by imposing $a_{i,i}=0$, $i=1,\dots,N$ also for $N>2$.
The equality of the densities in requirement 6. leads to $f_i^++f_i^- = f_j^++f_j^-$, $i,j= 1,\dots,N$.
Since we can express the $f^+_i$ in terms of $f_i^-$, we obtain after some direct computations the following constraint
\begin{align*}
 \sum_{j=1}^{N}a_{i,j} = 1 \qquad \forall j = 1,\dots,N\ .
\end{align*}
Applying all the above constraints in the case of a three-way junction $N=3$, 
only one free parameter $\alpha\in [0,1]$ is left in the entries of the matrix $A$
\begin{align*}
 A
 =\left(\begin{array}{ccc}
  0&\alpha&1-\alpha \\
  1-\alpha &0&\alpha\\
  \alpha&1-\alpha&0
  \end{array}\right)
  \ . 
\end{align*}
The only choice in which all edges are treated equally is $\alpha=\frac12$.
This leads to the very simple set of coupling conditions
\begin{equation}\label{eq:kinetic_cc3}
\left[
\begin{array}{c}
f_1^+\\
f_2^+\\
f_3^+\\
\end{array}
\right] = \left[
\begin{array}{ccc}
0&1/2&1/2\\
1/2&0&1/2\\
1/2&1/2&0\\
\end{array}\right]\left[
\begin{array}{c}
f_1^-\\
f_2^-\\
f_3^-\\
\end{array}
\right]\ .
\end{equation}
In the cases $N>3$ some more freedom in the choices of the $a_{i,j}$ is given,
but the values $a_{i,j}=\frac{1}{N-1}$ $i\neq j$ and $a_{i,i}=0$ remain an admissible choice.

\subsection{Coupling condition for half-moment closure}
Now we want to derive coupling conditions for the half moment system \eqref{equ:sc half-moment system} from the kinetic model.
Analogous to \eqref{eq:halfmoment_average} we define the following averaged quantity on each edge 
$\rho^+_i,\rho^-_i,q^+_i,q^-_i$ for $i = 1,\dots,N$.

Since the kinetic coupling conditions \eqref{eq:kinetic_cc3} are linear and independent of $v$
we obtain directly the coupling condition for half moment model \eqref{equ:sc half-moment system} as
\begin{align}\label{equ:coup_half_moment}
\left\{
\begin{array}{lcl}
\left[
\begin{array}{c}
\rho_1^+\\
\rho_2^+\\
\rho_3^+\\
\end{array}
\right] &=& \quad \left[
\begin{array}{ccc}
0&1/2&1/2\\
1/2&0&1/2\\
1/2&1/2&0\\
\end{array}\right]\left[
\begin{array}{c}
\rho_1^-\\
\rho_2^-\\
\rho_3^-\\
\end{array}
\right]
\\[0.5cm]
\left[
\begin{array}{c}
q_1^+\\
q_2^+\\
q_3^+\\
\end{array}
\right] &=& - \left[
\begin{array}{ccc}
0&1/2&1/2\\
1/2&0&1/2\\
1/2&1/2&0\\
\end{array}\right]\left[
\begin{array}{c}
q_1^-\\
q_2^-\\
q_3^-\\
\end{array}
\right]\\
\end{array}
\right.
\end{align}
These are six equations for six outgoing characteristics, which is the correct number of coupling conditions.
Note that only the averaging in $v$ and not the exact structure of the closure \eqref{eq:halfmoment_closure} was used.
If the coupling conditions depend on $v$ the closure \eqref{eq:halfmoment_closure} might be relevant.

Recall that all properties of \eqref{eq:kinetic_cc3} are inherited by the above coupling conditions.
For example the total mass in the system is conserved since
$$\dis\sum_{i=1}^3q_i = \dis\sum_{i=1}^3(q_i^+ + q_i^-)=0\ .$$
But note that although the coupling conditions maintain the positivity of the densities,
this does not necessarily hold for the complete network, as it is not assured by the model on the edges \eqref{equ:sc half-moment system}.

\subsection{Coupling condition for Cattaneo (p1) equations}\label{sec:CattaneoCoupling}
In order to derive coupling conditions for the macroscopic quantities $\rho$ and $q$ we can not simply average the equations \eqref{eq:kinetic_cc3} on $[-1,1]$, since the information is split for positive and negative values of $v$.
But from the closure \eqref{eq:halfmoment_closure} we can deduce the following expressions for the half moments
\begin{align*}
f^+_i=f(v) =\frac{1}{2}\rho_i + \frac{3}{2}v\epsilon q_i\ ,
&&
f^-_i=f(-v) =\frac{1}{2}\rho_i - \frac{3}{2}v\epsilon q_i,\quad v\in[0,1]\ ,\ i=1,\dots,N\ .
\end{align*}
Inserting these into \eqref{eq:kinetic_cc3} we obtain for the case $N=3$
\begin{equation}\label{eq:cc_P1}
 \left[
\begin{array}{ccc}
2&-1&-1\\
-1&2&-1\\
-1&-1&2\\
\end{array}\right] \left[
\begin{array}{c}
\rho_1\\
\rho_2\\
\rho_3\\
\end{array}
\right] + 
\epsilon \frac{3}{2}
\left[
\begin{array}{ccc}
2&1&1\\
1&2&1\\
1&1&2\\
\end{array}\right]\left[
\begin{array}{c}
q_1\\
q_2\\
q_3\\
\end{array}
\right]
=0\ .
\end{equation}
Note that in \eqref{equ:p1} we have only one characteristic moving to the right.
Thus for a node connecting three outgoing edges, we have to provide exactly three coupling conditions.
The matrix in front of $\rho$ only has rank $2$ whereas the matrix for $q$ has full rank.
How bad condition numbers for small $\epsilon$ can be avoided is indicated in section \ref{sec:CC_KellerSegel}.


As before all properties of the kinetic coupling conditions transfer to \eqref{eq:cc_P1}, 
e.g. the classical formulation of the conservation of mass is obtained by summing all three equations.

\subsubsection*{Alternative coupling conditions for the P1 model}
In \cite{2014arXiv1411.6109R} a type of coupling conditions for \eqref{equ:p1} has been proposed.
For only outgoing edges and the notation used in the present models these are of the form 
\begin{align*}
 \frac{1}{\sqrt{3}\epsilon}q_i = \sum_{j=1}^N\alpha_{i,j}\left(\rho_j-\rho_i\right)\qquad,\ i = 1,\dots,N 
\end{align*}
for suitable parameters $\alpha_{i,j}$.
The coupling conditions \eqref{eq:cc_P1} can be transformed in the above form when choosing for the parameters $\alpha_{i,j}$ the following values $\alpha_{i,j}=\frac{2}{3\sqrt{3}\epsilon^2} \geq 0$.
This choice matches all conditions imposed in \cite{2014arXiv1411.6109R}, such as 
$\sum_{i=1}^N(\alpha_{i,j}-\alpha_{j,i})=0$ $j=1,\dots,N$ and $\alpha_{i,j}\geq 0$ $i,j=1,\dots,N$.

Another possible choice of coupling conditions for \eqref{equ:p1} can be motivated similarly to \cite{GasNetworks} by assuming the continuity of the density across the node 
\begin{align}\label{eq:cc_P1_equalrho}
 \sum_{i=1}^N q_i=0\ ,&& \rho_i=\rho_j\ ,\ i,j=1,\dots,N\ ,\ i\neq j\ .
\end{align}
These conditions differ from \eqref{eq:cc_P1} and will be investigated in the numerical examples in section \ref{sec:numerical_tests}.

\subsection{Coupling condition for the Keller-Segel equations}\label{sec:CC_KellerSegel}
For the Keller-Segel model \eqref{equ:d_limit} we have to verify 
that the proposed coupling conditions converge for $\epsilon\rightarrow 0$ to those introduced in \cite{MR3149314}. 
These conditions are 
\begin{align}\label{eq:cc_KellerSegel}
 \sum_{i=1}^N q_i&=0
 \\ \nonumber
 \rho_i&=\rho_j\qquad ,\ i,j=1,\dots,N\ ,\ i\neq j\ ,
\end{align}
where $q_i = \dis\frac{1}{3\lambda}(\partial_x\rho_i) - \dis\frac{\alpha}{3\lambda}\dis\overline{\partial_x m}\rho_i$ as in \eqref{equ:d_limit}.

As depicted in figure \ref{fig:Model_Limits} there are three possible ways to reach \eqref{eq:cc_KellerSegel} starting from \eqref{eq:kinetic_cc3}.
The first one is a direct limit from the kinetic model to the diffusion limit.
Since in deriving the coupling conditions for the half moment model \eqref{equ:coup_half_moment} only the averaging was used and no closure assumption, the way from \eqref{eq:kinetic_cc3} to \eqref{eq:cc_KellerSegel} directly passes \eqref{equ:coup_half_moment}.
Thus we can start with the coupling conditions for the half moment model. 
By inserting the definitions \eqref{eq:halfmoment_defrhoqrhohatqhat} into \eqref{equ:coup_half_moment} we obtain 
\begin{align*}
\left\{
\begin{array}{lcl}
\left[\begin{array}{ccc}
-2&1&1\\
1&-2&1\\
1&1&-2\\
\end{array}\right]
\left[
\begin{array}{c}
\rho_1\\
\rho_2\\
\rho_3\\
\end{array}
\right]
- \epsilon
\left[
\begin{array}{ccc}
2&1&1\\
1&2&1\\
1&1&2\\
\end{array}\right]
\left[
\begin{array}{c}
\hat\rho_1\\
\hat\rho_2\\
\hat\rho_3\\
\end{array}
\right]
&=& 0
\\[0.5cm]
\left[\begin{array}{ccc}
-2&1&1\\
1&-2&1\\
1&1&-2\\
\end{array}\right]
\left[
\begin{array}{c}
\hat q_1\\
\hat q_2\\
\hat q_3\\
\end{array}
\right]
- \epsilon
\left[\begin{array}{ccc}
2&1&1\\
1&2&1\\
1&1&2\\
\end{array}\right]
\left[\begin{array}{c}
q_1\\
q_2\\
q_3\\
\end{array}\right]
&=& 0
\end{array}
\right.
\ .
\end{align*}
We introduce the shorter vector notation 
\begin{align*}
 R\ \rho -\epsilon S \hat\rho&=0\\
 R\ \hat q -\epsilon S q&=0\ ,
\end{align*}
with suitable choices of $R$ and $S$.
Note that the matrix $R$ has only rank $2$ whereas $S$ has full rank.
By setting $\epsilon=0$ we would loose two of the given six equations.
In order to find the correct representation for the limit process we transform the system into  
\begin{align*}
 \tilde R\ \rho -\epsilon \tilde S \hat\rho&=0\\
 \tilde R\ \hat q -\epsilon \tilde S q&=0\ ,
\end{align*}
where $\tilde R=T\ R$ and $\tilde S=T\ S$ with
\begin{align}\label{eq:matrix_T}
T=
\left[\begin{array}{ccc}
1&0&0\\
0&1&0\\
1&1&1\\
\end{array}\right]\ .
\end{align}
The matrix $\tilde R$ has now only zeros in the last row, such that we can eliminate the $\epsilon$ in this equation.
The equivalent coupling conditions read
\begin{align*}
\left\{
\begin{array}{lcl}
\left[\begin{array}{ccc}
-2&1&1\\
1&-2&1\\
0&0&0\\
\end{array}\right]
\left[
\begin{array}{c}
\rho_1\\
\rho_2\\
\rho_3\\
\end{array}
\right]
- 
\left[\begin{array}{ccc}
2\epsilon&\epsilon&\epsilon\\
\epsilon&2\epsilon&\epsilon\\
1&1&1\\
\end{array}\right]
\left[
\begin{array}{c}
\hat\rho_1\\
\hat\rho_2\\
\hat\rho_3\\
\end{array}
\right]
&=& 0
\\[0.5cm]
\left[\begin{array}{ccc}
-2&1&1\\
1&-2&1\\
0&0&0\\
\end{array}\right]
\left[
\begin{array}{c}
\hat q_1\\
\hat q_2\\
\hat q_3\\
\end{array}
\right]
- 
\left[\begin{array}{ccc}
2\epsilon&\epsilon&\epsilon\\
\epsilon&2\epsilon&\epsilon\\
1&1&1\\
\end{array}\right]
\left[\begin{array}{c}
q_1\\
q_2\\
q_3\\
\end{array}\right]
&=& 0
\ .
\end{array}
\right.
\end{align*}
By setting $\epsilon=0$ we now obtain 
\begin{align*}
 \rho_1=
 \rho_2&=
 \rho_3\ ,
 &&&
 \sum_{i=1}^3\hat \rho_i&=0\ ,\\
 \hat q_1=
 \hat q_2&=
 \hat q_3\ ,
 &&&
 \sum_{i=1}^3q_i&=0\ ,
\end{align*}
which contains the desired equations \eqref{eq:cc_KellerSegel}.

Similarly we can do the limit starting from the Cattaneo model \eqref{eq:cc_P1}.
Again we first have to apply the transformation \eqref{eq:matrix_T} before obtaining directly \eqref{eq:cc_KellerSegel} when setting $\epsilon =0$.

\section{Numerical schemes on an interval}\label{sec:num_interval}
For approximating the solutions of the models discussed in section \ref{sec:models} many different numerical schemes can be used. 
Since in all models the equation describing the chemoattractant is only coupled via source terms to the dynamics of the cells, we can handle both parts separately using a splitting
algorithm.
The equation for the chemoattractant $m$ is always solved with a forward central difference scheme as in \cite{MR3149314}.

Therefore, we address in the following only the discretization of the models for the movements of the cell density.
Except from the Keller-Segel model \eqref{equ:d_limit}, the schemes are based on the relaxation methods proposed in \cite{MR1322811} and \cite{MR1655853}.
For more sophisticated relaxation schemes for kinetic and Cattaneo model we refer to \cite{Gosse,Gosse2,refId0,2014arXiv1411.6109R}.
The main advantage of this scheme is that they are asymptotic preserving, i.e. they provide good approximations for any value of $\epsilon$ without severe time step restrictions. 

Each system is split into a relaxation part and a transport part.
The  system of ODEs resulting from the relaxation part is solved  implicitly with the backward Euler method in order to avoid any time step restriction depending on $\epsilon$.
The advection part considering the terms on the left hand sides is  solved  with the explicit upwind method.
Here a mild restriction on the time step due to the CFL condition arises.

\subsection{The kinetic model}
For the discretization of the   kinetic equation  \eqref{equ:g_kinetic} we consider the even- and odd-parities \eqref{cond}
 for positive velocities $(v\geq 0)$ and introduce additional terms following \cite{MR1655853}.
 This leads to the system
 \begin{align*}
 \partial_t r + v\partial_x j&=-\dis\frac{1}{\epsilon^2}\left(\lambda r-\frac{\lambda}{2}\rho\right) \\
 \partial_t j +  \phi v \partial_x r&=-\dis\frac{1}{\epsilon^2}\left(\lambda j -\frac{1}{2}\alpha v 
\overline{\partial_x m} \rho  +
(1- \epsilon^2 \phi)  v \partial_x r
 \right) 
 \end{align*}
 with $0 \le \phi \le  \frac{1}{ \epsilon^2}$ and $\epsilon \le \frac{\lambda}{\alpha}$.
 For example we may choose
 $\phi= \frac{\alpha^2}{\lambda^2}$.
 This system we discretize in $v$ direction with $N_v$ points, such that we obtain a system of $N_v$ equations for the relaxation scheme. 

\subsection{The $P_1$-model}
 Similar to the above model we introduce for \eqref{equ:p1} with
and $\epsilon \le \frac{\lambda}{\alpha}$ the speed
$0 \le \phi\le \frac{1}{3 \epsilon^2}$  and obtain
\begin{align*}
 \partial_t\rho + \partial_x q
 &=0 \\
 \partial_t q +\phi\partial_x\rho 
 &= -\dis\frac{1}{\epsilon^2}\left( \lambda q +\left(\frac{1}{3}-\epsilon^2\phi\right)\partial_x\rho - \frac{\alpha}{3}
 \overline{\partial_x m}
 \rho\right)
 \ .
 \end{align*}
 We choose $\phi = \frac{\alpha^2}{3 \lambda^2}$.
 \subsection{The half-moment model}
 We rearrange \eqref{equ:sc half-moment system_hat}
 and introduce new terms involving the speed $0 \le \phi \le  \frac{1}{6 \epsilon^2}$
 for  $\epsilon \le \frac{\lambda}{\alpha}$. We obtain
 \begin{align}\label{equ:relaxation_sys_U}
 \begin{array}{rl}
\partial_t\rho + \partial_xq
&= 0\\[0.3cm]
\partial_tq + \phi\partial_x\left(-\rho + 6\hat{q}\right) 
&=-\dis\frac{1}{\epsilon^2}\left(\lambda q - \dis\frac{\alpha}{3}\overline{\partial_x m}\rho + \left(\dis\frac{1}{6}-\epsilon^2\phi\right)\partial_x\left(-\rho +6\hat{q}\right)\right)\\[0.3cm]
\partial_t\hat{\rho} + \phi\partial_x\hat{q} 
&=  -\dis\frac{1}{\epsilon^2}\left(\lambda\hat{\rho} -\dis\frac{1}{2}\alpha \overline{\partial_x m}\rho+\left(1-\epsilon^2\phi\right)\partial_x\hat{q}\right)\\[0.3cm]
\partial_t\hat{q} + \dis\frac{1}{6}\partial_x\left(-\hat{\rho} + 6q\right) 
&= -\dis\frac{1}{\epsilon^2}\left(\hat{q}-\dis\frac{\rho}{2}\right).
 \end{array}
\end{align}
In this case we choose $\phi = \frac{\alpha^2}{6 \lambda^2}$.
This system can be directly solved with the relaxation scheme.

\subsection{Keller-Segel equation}
The Keller-Segel \eqref{equ:d_limit} equation we discretize as in \cite{MR3149314} by a central difference approximation in space and the forward Euler method in time.
Note that the discretizations of the above models yield a consistent and stable approximation of 
\begin{equation*}
 \partial _t \rho -\partial_x\left(\dis\frac{1}{3\lambda}\partial_x\rho-\dis\frac{\alpha}{3\lambda}\dis\frac{\partial_x m}{\sqrt{1+\vert\partial_xm\vert^2}}\rho\right) = 0
\end{equation*}
in the case $\epsilon=0$. 
These discretizations are different from the one used for the Keller-Segel model, in particular,  they  have a wider numerical stencil. 
 
\section{Numerical methods for chemotaxis models on networks}\label{sec:numerical_coupling}
In this section we discuss how the above numerical methods can be extended to models for chemotaxis on networks.
This includes a suitable discretization of the coupling conditions introduced in section \ref{sec:CouplingConditions}.
As the equation for the chemoattractant is splitted from the remaining equations also the coupling procedure can be regarded separately.
A numerical procedure for the coupling of $m$ has been introduced in \cite{MR3149314} and is also used in the test cases of the present paper.
As noted in section \ref{sec:CouplingConditions} for the equations describing the cell movement only the advective component contribute to the coupling procedure. 
In the schemes of section \ref{sec:num_interval} the advection step is clearly separated from the relaxation step. 
Therefore we address the numerical coupling procedures only for their first step, since in the relaxation step no coupling is needed.

As in section \ref{sec:CouplingConditions} mainly the case of a junction coupling three outgoing edges will be considered, 
but all procedures easily extend to arbitrary junctions.
Based on these tools networks of any shape can be handled, just by connecting the edges with nodes accordingly. 

\subsection{Numerical coupling of the kinetic equation}
We consider the advection step for each discrete velocity $v_k$, $k = 1,\dots,N_v$
\begin{align}\label{equ:convection_f}
\begin{array}{rl}
 \partial_t r_k + v_k(\partial_x j_k)&=0 \\
 \partial_t j_k + \phi v_k (\partial_x r_k)&=0 \ .
\end{array}
\end{align}
Inserting \eqref{eq:evenodd} into \eqref{eq:kinetic_cc3} yields for every $v_k$ the coupling conditions
\begin{equation}\label{kinetic_cc_rj}
\left[
\begin{array}{c}
r_1^i + \epsilon j_1^i\\
r_2^i + \epsilon j_2^i\\
r_3^i + \epsilon j_3^i\\
\end{array}
\right] 
= \left[
\begin{array}{ccc}
0&1/2&1/2\\
1/2&0&1/2\\
1/2&1/2&0\\
\end{array}\right]\left[
\begin{array}{c}
r_1^i - \epsilon j_1^i\\
r_2^i - \epsilon j_2^i\\
r_3^i - \epsilon j_3^i\\
\end{array}
\right]
\end{equation}
at the node.
For a correct treatment of these boundary values we decompose \eqref{equ:convection_f} into characteristic variables, which are related to the eigenvectors of the transport matrix in  \eqref{equ:convection_f}. 
Since the characteristic variables associated with negative wave speeds exit the domain,
we can use \eqref{kinetic_cc_rj} to determine the values of the characteristic variables for the positive eigenvalues $v_k+\sqrt{\phi}$.
For each connected edge and each $v_k$ we have one wave with positive speed, which yields three unknowns in three equations \eqref{kinetic_cc_rj}.
Once the values of the characteristic variables are determined, we can compute the values of $r_k$ and $j_k$ at the junction.

\subsection{Numerical coupling of the half moment equations}
For the half moment model we apply a similar technique as for the kinetic equation.
Note that the half moment model \eqref{equ:sc half-moment system} is already decomposed into forward and backward going components.
But the relaxation scheme solves \eqref{equ:relaxation_sys_U} and since the wave speeds are modified by the relaxation, we can not use the inverse of the transformations \eqref{eq:halfmoment_defrhoqrhohatqhat} to separate in and outgoing waves.

Thus we express the coupling conditions \eqref{equ:coup_half_moment} in the variables of \eqref{equ:sc half-moment system_hat} 
\begin{align*}
\begin{array}{lcl}
\left[
\begin{array}{c}
\rho_1+\epsilon \hat\rho_1\\
\rho_2+\epsilon \hat\rho_2\\
\rho_3+\epsilon \hat\rho_3
\end{array}
\right] &=& \quad \left[
\begin{array}{ccc}
0&1/2&1/2\\
1/2&0&1/2\\
1/2&1/2&0\\
\end{array}\right]\left[
\begin{array}{c}
\rho_1-\epsilon \hat\rho_1\\
\rho_2-\epsilon \hat\rho_2\\
\rho_3-\epsilon \hat\rho_3
\end{array}
\right]
\\[0.5cm]
\left[
\begin{array}{c}
\epsilon q_1+ \hat q_1\\
\epsilon q_2+ \hat q_2\\
\epsilon q_3+ \hat q_3
\end{array}
\right] &=& - \left[
\begin{array}{ccc}
0&1/2&1/2\\
1/2&0&1/2\\
1/2&1/2&0\\
\end{array}\right]\left[
\begin{array}{c}
\epsilon q_1- \hat q_1\\
\epsilon q_2- \hat q_2\\
\epsilon q_3- \hat q_3
\end{array}
\right]\ .
\end{array}
\end{align*}
Now we again consider the characteristic decomposition of \eqref{equ:relaxation_sys_U} and use the above equations to determine the strengths of the forward going waves.
Since on each edge there are two waves with positive speed, we need exactly six equations for three edges, as provided.

\subsection{Numerical coupling of the Cattaneo equations}
Independent of the choice of coupling conditions for the Cattaneo equations in section \ref{sec:CattaneoCoupling},
we can apply the same numerical procedure. 
As above we decompose the advective part of the relaxed system into characteristic variables and fix the values along the outgoing wave.
The strength of the ingoing waves is easily determined by the coupling conditions.

\section{Numerical tests}\label{sec:numerical_tests}
In this section we investigate the proposed models in several numerical test cases.
A main focus will be given to the behavior of the models for different values of $\epsilon$.

In all the considered examples we will use the following values for the parameters
$\lambda =\alpha =1$, $D = 1$, $\gamma_{\rho}=1$ and $\gamma_m=0.1$. 
The spatial resolution is $\Delta x = 0.02$ and the time step is chosen according to the CFL condition.
At end points where no coupling conditions are imposed zero Neumann boundary conditions are applied.
For the kinetic model we discretize the velocity space $V=[-1,1]$ with $N_v=50$ cells.
In the comparisons the model proposed in \cite{2014arXiv1411.6109R} is shown with the parameters $\alpha_{i,j}=1$.

\subsection{Numerical solutions on an interval}
\begin{figure}[ht!]
\externaltikz{figure_interval_1}{
\centering
\ifthenelse{\isundefined{\pfad}}{
  \def\pfad{Photos_tikz/numerical_tests/Data/interval/}
}{
}

\begin{tikzpicture}
\begin{groupplot}[group style={group size=2 by 1, horizontal sep = 2cm,  vertical sep = 2cm},
width=10cm,
height=4cm,
scale only axis,
xmin = 0,
xmax = 2,
ymin = 0,
]
\nextgroupplot[
title = {$\epsilon = 1$},
xlabel = {$x$},
ylabel = {$\rho$},
legend cell align=left,
]
\addplot+ [mark repeat = 500, mark phase = 215, smooth,line width=2pt] table[x = x,y = r_kinetic] {\pfad line_e1.txt};
\addlegendentry{Kinetic};

\addplot+ [mark repeat = 500, mark phase = 189,line width=2pt,smooth] table[x = x,y = r_half_M] {\pfad line_e1.txt};
\addlegendentry{Half-moment};

\addplot+ [mark repeat = 500, mark phase = 173,line width=2pt,smooth] table[x = x,y = r_cattaneo] {\pfad line_e1.txt};
\addlegendentry{Cattaneo};

\addplot+ [mark repeat = 500, mark phase = 180,line width=2pt,smooth] table[x = x,y = r_keller] {\pfad line_e1.txt};
\addlegendentry{Keller-Segel};

\addplot+ [mark repeat = 500, mark phase = 190,line width=2pt] table[x = x,y = r_initial] {\pfad line_e1.txt};
\addlegendentry{Initial condition};
\end{groupplot}
\end{tikzpicture}
}\\
\externaltikz{figure_interval_05}{
\ifthenelse{\isundefined{\pfad}}{
  \def\pfad{Photos_tikz/numerical_tests/Data/interval/}
}{
}

\begin{tikzpicture}
\begin{groupplot}[group style={group size=2 by 1, horizontal sep = 2cm,  vertical sep = 2cm},
width=10cm,
height=4cm,
scale only axis,
xmin = 0,
xmax = 2,
ymin = 0,
]
\nextgroupplot[
title = {$\epsilon = 0.5$},
xlabel = {$x$},
ylabel = {$\rho$},
legend cell align=left,
]
\addplot+ [mark repeat = 500, mark phase = 215, smooth,line width=2pt] table[x = x,y = r_kinetic] {\pfad line_e05.txt};
\addlegendentry{Kinetic};

\addplot+ [mark repeat = 500, mark phase = 180,line width=2pt,smooth] table[x = x,y = r_half_M] {\pfad line_e05.txt};
\addlegendentry{Half-moment};

\addplot+ [mark repeat = 500, mark phase = 173,line width=2pt,smooth] table[x = x,y = r_cattaneo] {\pfad line_e05.txt};
\addlegendentry{Cattaneo};

\addplot+ [mark repeat = 500, mark phase = 180,line width=2pt,smooth] table[x = x,y = r_keller] {\pfad line_e05.txt};
\addlegendentry{Keller-Segel};
\addplot+ [mark repeat = 500, mark phase = 190,line width=2pt] table[x = x,y = r_initial] {\pfad line_e05.txt};
\addlegendentry{Initial condition};
\end{groupplot}
\end{tikzpicture}
}
\caption{Numerical solutions of the four models on an interval at time $t=0.2$ with $\epsilon = 1$ (top) and $\epsilon = 0.5$ (bottom)}
\label{fig:interval_e1_e05}
\end{figure}
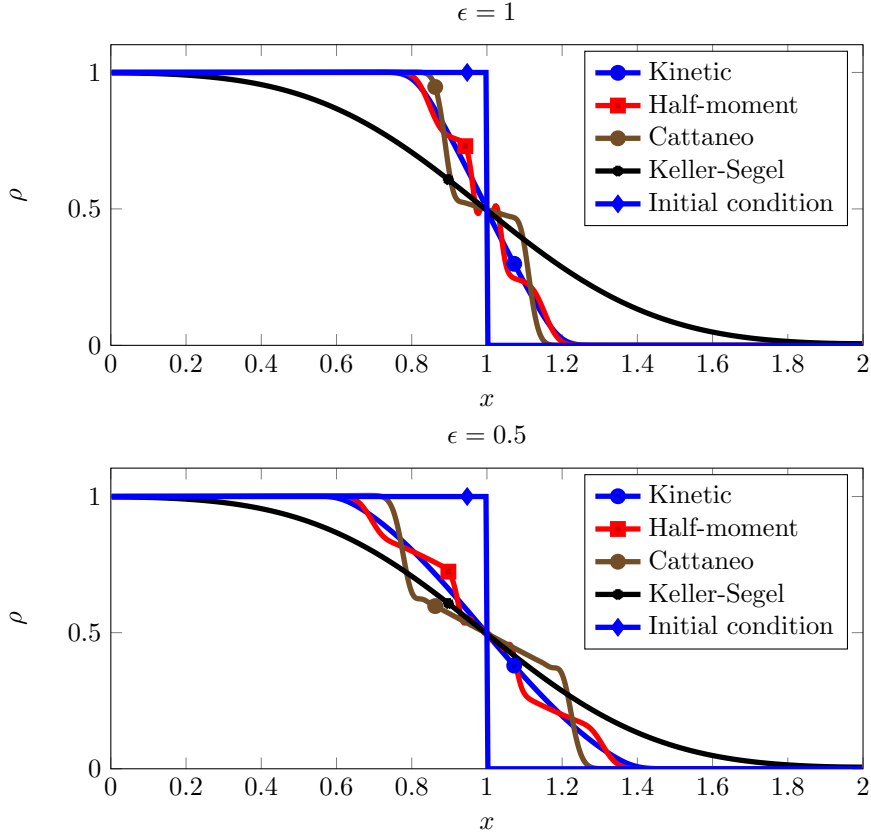
\begin{figure}[ht!]
\externaltikz{figure_interval_01}{
\centering
\ifthenelse{\isundefined{\pfad}}{
  \def\pfad{Photos_tikz/numerical_tests/Data/interval/}
}{
}

\begin{tikzpicture}
\begin{groupplot}[group style={group size=2 by 1, horizontal sep = 2cm,  vertical sep = 2cm},
width=10cm,
height=4cm,
scale only axis,
xmin = 0,
xmax = 2,
ymin = 0,
]
\nextgroupplot[
title = {$\epsilon = 0.1$},
xlabel = {$x$},
ylabel = {$\rho$},
legend cell align=left,
]
\addplot+ [mark repeat = 500, mark phase = 215, smooth,line width=2pt] table[x = x,y = r_kinetic] {\pfad line_e01.txt};
\addlegendentry{Kinetic};

\addplot+ [mark repeat = 500, mark phase = 189,line width=2pt,smooth] table[x = x,y = r_half_M] {\pfad line_e01.txt};
\addlegendentry{Half-moment};

\addplot+ [mark repeat = 500, mark phase = 173,line width=2pt,smooth] table[x = x,y = r_cattaneo] {\pfad line_e01.txt};
\addlegendentry{Cattaneo};

\addplot+ [mark repeat = 500, mark phase = 180,line width=2pt,smooth] table[x = x,y = r_keller] {\pfad line_e01.txt};
\addlegendentry{Keller-Segel};
\addplot+ [mark repeat = 500, mark phase = 190,line width=2pt] table[x = x,y = r_initial] {\pfad line_e01.txt};
\addlegendentry{Initial condition};
\end{groupplot}
\end{tikzpicture}
}\\
\externaltikz{figure_interval_0}{
\ifthenelse{\isundefined{\pfad}}{
  \def\pfad{Photos_tikz/numerical_tests/Data/interval/}
}{
}

\begin{tikzpicture}
\begin{groupplot}[group style={group size=2 by 1, horizontal sep = 2cm,  vertical sep = 2cm},
width=10cm,
height=4cm,
scale only axis,
xmin = 0,
xmax = 2,
ymin = 0,
]
\nextgroupplot[
title = {$\epsilon \approx 0$},
xlabel = {$x$},
ylabel = {$\rho$},
legend cell align=left,
]
\addplot+ [mark repeat = 500, mark phase = 215, smooth,line width=2pt] table[x = x,y = r_kinetic] {\pfad line_e0.txt};
\addlegendentry{Kinetic};

\addplot+ [mark repeat = 500, mark phase = 189,line width=2pt,smooth] table[x = x,y = r_half_M] {\pfad line_e0.txt};
\addlegendentry{Half-moment};

\addplot+ [mark repeat = 500, mark phase = 173,line width=2pt,smooth] table[x = x,y = r_cattaneo] {\pfad line_e0.txt};
\addlegendentry{Cattaneo};

\addplot+ [mark repeat = 500, mark phase = 180,line width=2pt,smooth] table[x = x,y = r_keller] {\pfad line_e0.txt};
\addlegendentry{Keller-Segel};

\addplot+ [mark repeat = 500, mark phase = 190,line width=2pt] table[x = x,y = r_initial] {\pfad line_e0.txt};
\addlegendentry{Initial condition};
\end{groupplot}
\end{tikzpicture}
}
\caption{Numerical solutions of the four models on an interval at time $t=0.2$ with $\epsilon = 0.1$ (top) and $\epsilon = 10^{-6}$ (bottom)}
\label{fig:interval_e01_e0}
\end{figure}
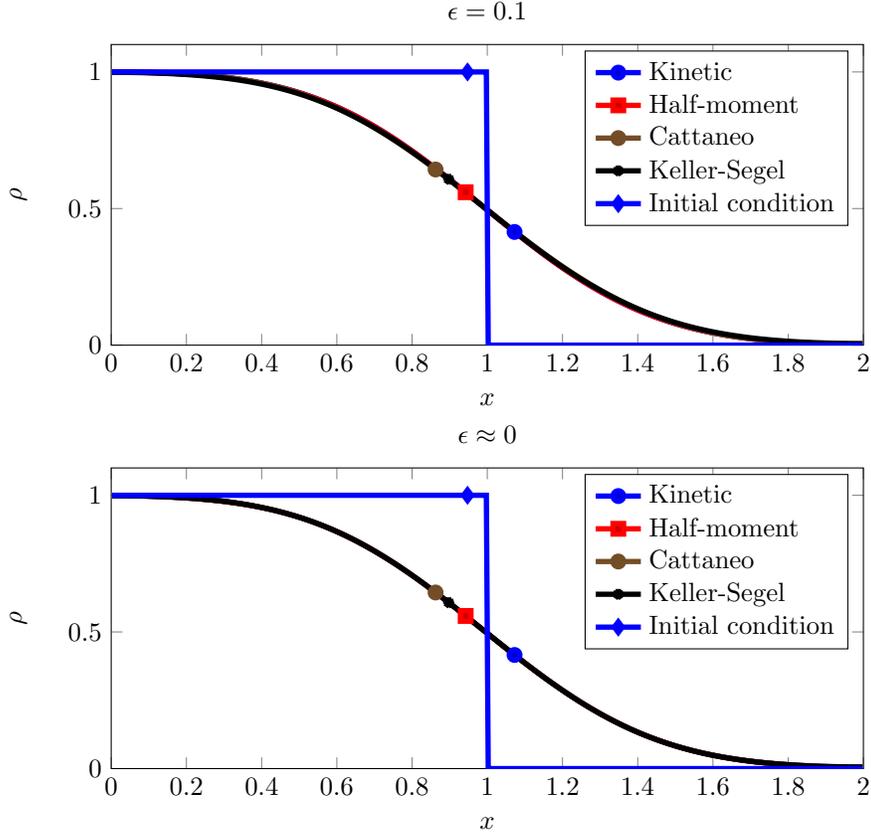
In these test cases we study the elementary behaviour of the different models, especially their convergence as $\epsilon \rightarrow 0$.
Consider the interval $[0,2]$ and $\Delta x=0.005$.
As initial conditions a Riemann problem is given with
\begin{equation*}
\rho(x,0)=
\left\{
\begin{array}{lr}
1&\text{ if } 0\leq x\leq 1 \\
0&\text{ if } 1 \leq x\leq 2 
\end{array}
\right.
\end{equation*}
such that for kinetic equation we obtain
\begin{equation*}
f(x,v,0)=F(v)\rho(x,0)=\dis\frac{1}{2}\rho(x,0)\ .
\end{equation*}
The corresponding initial values for the macroscopic equations can be derived thereof.
At $t=0$ no chemoattractant $m(x,0)=0$ is present.

In the Figures \ref{fig:interval_e1_e05}-\ref{fig:interval_e01_e0} the densities at $t=0.2$ for the values $\epsilon=1,\ 0.5,\ 0.1,\ 10^{-6}$ are shown.
For the kinetic model we observe that the speed of the movement strongly depends on the value of $\epsilon$,
the smaller $\epsilon$ is the faster the waves propagate.
As expected the closest approximation to the kinetic equations is the half moment model.
For  $\epsilon=1$ and  $\epsilon=0.5$ one can clearly observe the four waves generated by the advective part.
In the Cattaneo model only two waves can be used to approximate the kinetic wave.
The Keller-Segel equation is evolving too fast for $\epsilon$ large. 
For small $\epsilon$ all models converge to the pattern of the Keller-Segel equation.

\subsection{Numerical solutions on a tripod network}
This test case focuses on the coupling conditions of the different models proposed in section \ref{sec:CouplingConditions}.
Therefore we study a junction connecting three outgoing edges, as depicted in Figure \ref{fig:out_tripod}.
On each of the edges we consider the interval $[0,1]$.
As initial conditions we choose 
\begin{align*}
\rho_1(x,0) = 1\ ,&&
\rho_2(x,0) = 2\ ,&&
\rho_3(x,0) = 3\ ,
\end{align*}
which yields the following values for kinetic equation
\begin{equation*}
f_i(x,v,0)=F(v)\rho_i(x,0)=\dis\frac{1}{2}\rho_i(x,0)\qquad i =1,2,3\ .
\end{equation*}
All other quantities are initially zero.

\begin{figure}[ht!]
\centering
\ifthenelse{\isundefined{\pfad}}{
  \def\pfad{./Photos_tikz/numerical_tests/Data/tripod_net/}
}{
}

\centering
	\begin{tikzpicture}[baseline]
		\begin{axis}
		[	width=6cm,
			height=7cm,
			scale only axis,
			xmin = 0,
			xmax = 1,
			ymin = 0.9,
			ymax = 3.1,
			title = {1st edge},
			xlabel = {$x$},
			ylabel = {$\rho$},
			name=Axis1,
			legend cell align=left,
			legend style={at={(0,-0.25)},
						  anchor=north west,
						  },			  
							]

			\addplot+ [mark repeat = 50, mark phase = 47, thick,smooth] table[x = x,y = r_kinetic] {\pfad edge_1st1.txt};
			\addlegendentry{Kinetic};

			\addplot+ [mark repeat = 49, mark phase = 46, thick,smooth] table[x = x,y = r_half] {\pfad edge_1st1.txt};
			\addlegendentry{Half-moment};
			
			\addplot+ [mark repeat = 49, mark phase = 44, thick,smooth] table[x = x,y = r_p1] {\pfad edge_1st1.txt};
			\addlegendentry{p1};
			\addplot+ [mark repeat = 49, mark phase = 47, thick,smooth] table[x = x,y = r_fag] {\pfad edge_1st1.txt};
			\addlegendentry{p1\cite{2014arXiv1411.6109R}};
			
			\addplot+ [mark repeat = 49, mark phase = 47, thick,smooth] table[x = x,y = r_con] {\pfad edge_1st1.txt};
			\addlegendentry{p1$(\rho_i=\rho_j)$};
			\addplot+ [mark repeat = 50, mark phase = 25, thick,smooth] table[x = x,y = r_keller] {\pfad edge_1st1.txt};
			\addlegendentry{Keller-Segel};
			\addplot+ [mark repeat = 49, mark phase = 44, thick,smooth] table[x = x,y = r_in1] {\pfad edge_1st1.txt};
			\addlegendentry{Initial condition};

		\end{axis}

		\begin{axis}[
			width=4cm,
			height=7cm,
			scale only axis,
			xmin = 0,
			xmax = 1,
			ymin = 0.9,
			ymax = 3.1,
			xshift=0.8cm,
			at={(Axis1.east)},
			anchor=west,
			name=Axis2,					
			title = {3rd edge},
			xlabel = $x$,
			]
							]

			\addplot+ [mark repeat = 60, mark phase = 5, thick,smooth] table[x = x,y = r_kinetic] {\pfad edge_3rd1.txt};

			\addplot+ [mark repeat = 60, mark phase = 5, thick,smooth] table[x = x,y = r_half] {\pfad edge_3rd1.txt};
			\addplot+ [mark repeat = 60, mark phase = 2, thick,smooth] table[x = x,y = r_p1] {\pfad edge_3rd1.txt};
			\addplot+ [mark repeat = 60, mark phase = 2, thick,smooth] table[x = x,y = r_fag] {\pfad edge_3rd1.txt};
			
			\addplot+ [mark repeat = 60, mark phase = 4, thick,smooth] table[x = x,y = r_con] {\pfad edge_3rd1.txt};
			
			\addplot+ [mark repeat = 60, mark phase = 20, thick,smooth] table[x = x,y = r_keller] {\pfad edge_3rd1.txt};
			\addplot+ [mark repeat = 60, mark phase = 5, thick,smooth] table[x = x,y = r_in3] {\pfad edge_3rd1.txt};
			
		\end{axis}
		
		\begin{axis}[
			width=4cm,
			height=4cm,
			scale only axis,
			xmin = 0,
			xmax = 1,
			yshift=-2cm,
			at={(Axis2.south)},
			anchor=north,
			name=Axis3,
			title = {2nd edge},
			xlabel = {$x$},
			ylabel = {$\rho$},
							]

			\addplot+ [mark repeat = 60, mark phase = 1, thick,smooth] table[x = x,y = r_kinetic] {\pfad edge_2nd1.txt};

			\addplot+ [mark repeat = 60, mark phase = 1, thick,smooth] table[x = x,y = r_half] {\pfad edge_2nd1.txt};
			\addplot+ [mark repeat = 60, mark phase = 1, thick,smooth] table[x = x,y = r_p1] {\pfad edge_2nd1.txt};
			\addplot+ [mark repeat = 60, mark phase = 1, thick,smooth] table[x = x,y = r_fag] {\pfad edge_2nd1.txt};
			
			\addplot+ [mark repeat = 60, mark phase = 1, thick,smooth] table[x = x,y = r_con] {\pfad edge_2nd1.txt};
			\addplot+ [mark repeat = 60, mark phase = 1, thick,smooth] table[x = x,y = r_keller] {\pfad edge_2nd1.txt};
			\addplot+ [mark repeat = 60, mark phase = 1, thick,smooth] table[x = x,y = r_in2] {\pfad edge_2nd1.txt};
			
		\end{axis}
		\node[anchor=center] at (Axis1.center) {\large$\epsilon=1$};

%
	\end{tikzpicture}%
\let\pfad\undefined
\caption{Numerical solutions on a tripod network at time $t=0.3$, $\Delta x=0.02$, $\epsilon = 1$}
\label{fig:tripod_e1}
\end{figure}
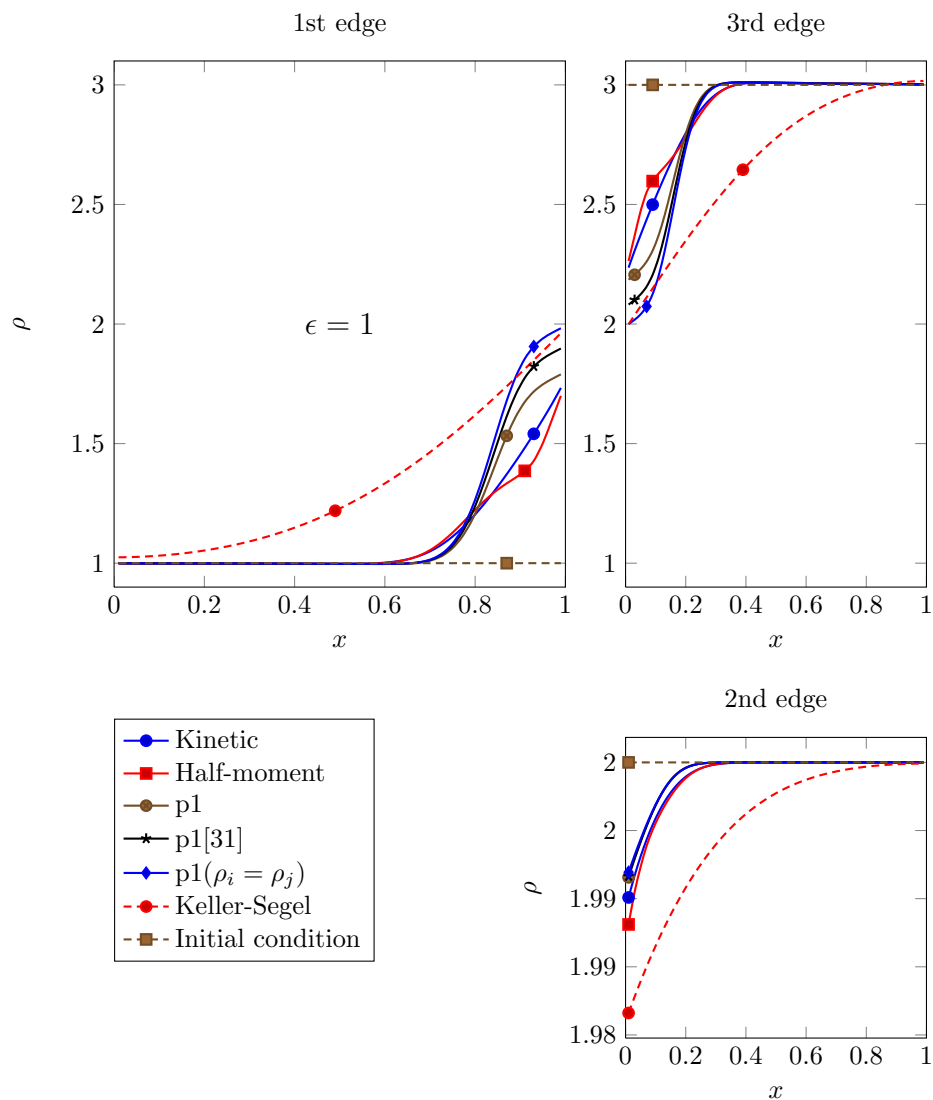
\begin{figure}[ht!]
\centering
\ifthenelse{\isundefined{\pfad}}{
  \def\pfad{Photos_tikz/numerical_tests/Data/tripod_net/}
}{
}

\newcommand{\subfile}{05}
\centering
	\begin{tikzpicture}[baseline]
		\begin{axis}
		[	width=6cm,
			height=7cm,
			scale only axis,
			xmin = 0,
			xmax = 1,
			ymin = 0.9,
			ymax = 3.1,
			title = {1st edge},
			xlabel = {$x$},
			ylabel = {$\rho$},
			name=Axis1,
			legend cell align=left,
			legend style={at={(0,-0.25)},
						  anchor=north west,
						  },			  
							]

			\addplot+ [mark repeat = 50, mark phase = 47, thick,smooth] table[x = x,y = r_kinetic] {\pfad edge_1st\subfile.txt};
			\addlegendentry{Kinetic};
			
			\addplot+ [mark repeat = 49, mark phase = 46, thick,smooth] table[x = x,y = r_half] {\pfad edge_1st\subfile.txt};
			\addlegendentry{Half-moment};
			
			\addplot+ [mark repeat = 49, mark phase = 44, thick,smooth] table[x = x,y = r_p1] {\pfad edge_1st\subfile.txt};
			\addlegendentry{p1};
			\addplot+ [mark repeat = 49, mark phase = 47, thick,smooth] table[x = x,y = r_fag] {\pfad edge_1st\subfile.txt};
			\addlegendentry{p1\cite{2014arXiv1411.6109R}};
			
			\addplot+ [mark repeat = 49, mark phase = 47, thick,smooth] table[x = x,y = r_con] {\pfad edge_1st\subfile.txt};
			\addlegendentry{p1$(\rho_i=\rho_j)$};
			
			\addplot+ [mark repeat = 50, mark phase = 25, thick,smooth] table[x = x,y = r_keller] {\pfad edge_1st\subfile.txt};
			\addlegendentry{Keller-Segel};
			
			\addplot+ [mark repeat = 49, mark phase = 44, thick,smooth] table[x = x,y = r_in1] {\pfad edge_1st\subfile.txt};
			\addlegendentry{Initial condition};

		\end{axis}

		\begin{axis}[
			width=4cm,
			height=7cm,
			scale only axis,
			xmin = 0,
			xmax = 1,
			ymin = 0.9,
			ymax = 3.1,
			xshift=0.8cm,
			at={(Axis1.east)},
			anchor=west,
			name=Axis2,					
			title = {3rd edge},
			xlabel = $x$,
			]
							]

			\addplot+ [mark repeat = 60, mark phase = 5, thick,smooth] table[x = x,y = r_kinetic] {\pfad edge_3rd\subfile.txt};
			
			\addplot+ [mark repeat = 60, mark phase = 5, thick,smooth] table[x = x,y = r_half] {\pfad edge_3rd\subfile.txt};
			\addplot+ [mark repeat = 60, mark phase = 2, thick,smooth] table[x = x,y = r_p1] {\pfad edge_3rd\subfile.txt};
			\addplot+ [mark repeat = 60, mark phase = 2, thick,smooth] table[x = x,y = r_fag] {\pfad edge_3rd\subfile.txt};
			
			\addplot+ [mark repeat = 60, mark phase = 4, thick,smooth] table[x = x,y = r_con] {\pfad edge_3rd\subfile.txt};
			\addplot+ [mark repeat = 60, mark phase = 20, thick,smooth] table[x = x,y = r_keller] {\pfad edge_3rd\subfile.txt};
			
			\addplot+ [mark repeat = 60, mark phase = 5, thick,smooth] table[x = x,y = r_in3] {\pfad edge_3rd\subfile.txt};
			
		\end{axis}
		
		\begin{axis}[
			width=4cm,
			height=4cm,
			scale only axis,
			xmin = 0,
			xmax = 1,
			yshift=-2cm,
			at={(Axis2.south)},
			anchor=north,
			name=Axis3,
			title = {2nd edge},
			xlabel = {$x$},
			ylabel = {$\rho$},
							]

			\addplot+ [mark repeat = 60, mark phase = 1, thick,smooth] table[x = x,y = r_kinetic] {\pfad edge_2nd\subfile.txt};

			\addplot+ [mark repeat = 60, mark phase = 1, thick,smooth] table[x = x,y = r_half] {\pfad edge_2nd\subfile.txt};
			\addplot+ [mark repeat = 60, mark phase = 1, thick,smooth] table[x = x,y = r_p1] {\pfad edge_2nd\subfile.txt};
			\addplot+ [mark repeat = 60, mark phase = 1, thick,smooth] table[x = x,y = r_fag] {\pfad edge_2nd\subfile.txt};
			
			\addplot+ [mark repeat = 60, mark phase = 1, thick,smooth] table[x = x,y = r_con] {\pfad edge_2nd\subfile.txt};
			\addplot+ [mark repeat = 60, mark phase = 1, thick,smooth] table[x = x,y = r_keller] {\pfad edge_2nd\subfile.txt};
			\addplot+ [mark repeat = 60, mark phase = 1, thick,smooth] table[x = x,y = r_in2] {\pfad edge_2nd\subfile.txt};
			
		\end{axis}
		\node[anchor=center] at (Axis1.center) {\large$\epsilon=0.5$};

%
	\end{tikzpicture}%
\let\pfad\undefined
\caption{Numerical solutions on a tripod network at time $t=0.3$, $\Delta x=0.02$, $\epsilon = 0.5$}
\label{fig:tripod_e05}
\end{figure}
\begin{figure}[ht!]
\centering
\ifthenelse{\isundefined{\pfad}}{
  \def\pfad{Photos_tikz/numerical_tests/Data/tripod_net/}
}{
}

\newcommand{\subfile}{01}
\centering
	\begin{tikzpicture}[baseline]
		\begin{axis}
		[	width=6cm,
			height=7cm,
			scale only axis,
			xmin = 0,
			xmax = 1,
			ymin = 0.9,
			ymax = 3.1,
			title = {1st edge},
			xlabel = {$x$},
			ylabel = {$\rho$},
			name=Axis1,
			legend cell align=left,
			legend style={at={(0,-0.25)},
						  anchor=north west,
						  },			  
							]

			\addplot+ [mark repeat = 50, mark phase = 47, thick,smooth] table[x = x,y = r_kinetic] {\pfad edge_1st\subfile.txt};
			\addlegendentry{Kinetic};
			
			\addplot+ [mark repeat = 49, mark phase = 46, thick,smooth] table[x = x,y = r_half] {\pfad edge_1st\subfile.txt};
			\addlegendentry{Half-moment};
			
			\addplot+ [mark repeat = 49, mark phase = 44, thick,smooth] table[x = x,y = r_p1] {\pfad edge_1st\subfile.txt};
			\addlegendentry{p1};
			\addplot+ [mark repeat = 49, mark phase = 47, thick,smooth] table[x = x,y = r_fag] {\pfad edge_1st\subfile.txt};
			\addlegendentry{p1 \cite{2014arXiv1411.6109R}};
			
			\addplot+ [mark repeat = 49, mark phase = 47, thick,smooth] table[x = x,y = r_con] {\pfad edge_1st\subfile.txt};
			\addlegendentry{p1 $(\rho_i=\rho_j)$};
			
			\addplot+ [mark repeat = 50, mark phase = 25, thick,smooth] table[x = x,y = r_keller] {\pfad edge_1st\subfile.txt};
			\addlegendentry{Keller-Segel};
			
			\addplot+ [mark repeat = 49, mark phase = 44, thick,smooth] table[x = x,y = r_in1] {\pfad edge_1st\subfile.txt};
			\addlegendentry{Initial condition};

		\end{axis}

		\begin{axis}[
			width=4cm,
			height=7cm,
			scale only axis,
			xmin = 0,
			xmax = 1,
			ymin = 0.9,
			ymax = 3.1,		
			xshift=0.8cm,
			at={(Axis1.east)},
			anchor=west,
			name=Axis2,					
			title = {3rd edge},
			xlabel = $x$,
			]
							]

			\addplot+ [mark repeat = 60, mark phase = 5, thick,smooth] table[x = x,y = r_kinetic] {\pfad edge_3rd\subfile.txt};

			\addplot+ [mark repeat = 60, mark phase = 5, thick,smooth] table[x = x,y = r_half] {\pfad edge_3rd\subfile.txt};
			\addplot+ [mark repeat = 60, mark phase = 2, thick,smooth] table[x = x,y = r_p1] {\pfad edge_3rd\subfile.txt};
			\addplot+ [mark repeat = 60, mark phase = 2, thick,smooth] table[x = x,y = r_fag] {\pfad edge_3rd\subfile.txt};
			
			\addplot+ [mark repeat = 60, mark phase = 4, thick,smooth] table[x = x,y = r_con] {\pfad edge_3rd\subfile.txt};
			\addplot+ [mark repeat = 60, mark phase = 20, thick,smooth] table[x = x,y = r_keller] {\pfad edge_3rd\subfile.txt};
			\addplot+ [mark repeat = 60, mark phase = 5, thick,smooth] table[x = x,y = r_in3] {\pfad edge_3rd\subfile.txt};
			
		\end{axis}
		
		\begin{axis}[
			width=4cm,
			height=4cm,
			scale only axis,
			xmin = 0,
			xmax = 1,
			yshift=-2cm,
			at={(Axis2.south)},
			anchor=north,
			name=Axis3,
			title = {2nd edge},
			xlabel = {$x$},
			ylabel = {$\rho$},
							]

			\addplot+ [mark repeat = 60, mark phase = 1, thick,smooth] table[x = x,y = r_kinetic] {\pfad edge_2nd\subfile.txt};

			\addplot+ [mark repeat = 60, mark phase = 1, thick,smooth] table[x = x,y = r_half] {\pfad edge_2nd\subfile.txt};
			\addplot+ [mark repeat = 60, mark phase = 1, thick,smooth] table[x = x,y = r_p1] {\pfad edge_2nd\subfile.txt};
			\addplot+ [mark repeat = 60, mark phase = 1, thick,smooth] table[x = x,y = r_fag] {\pfad edge_2nd\subfile.txt};
			
			\addplot+ [mark repeat = 60, mark phase = 1, thick,smooth] table[x = x,y = r_con] {\pfad edge_2nd\subfile.txt};
			\addplot+ [mark repeat = 60, mark phase = 1, thick,smooth] table[x = x,y = r_keller] {\pfad edge_2nd\subfile.txt};
			\addplot+ [mark repeat = 60, mark phase = 1, thick,smooth] table[x = x,y = r_in2] {\pfad edge_2nd\subfile.txt};
			
		\end{axis}
		
		\node[anchor=center] at (Axis1.center) {\large$\epsilon=0.1$};
	\end{tikzpicture}%
\let\pfad\undefined
\caption{Numerical solutions on a tripod network at time $t=0.3$, $\Delta x=0.02$, $\epsilon = 0.1$}
\label{fig:tripod_e01}
\end{figure}
\begin{figure}[ht!]
\centering
\ifthenelse{\isundefined{\pfad}}{
  \def\pfad{Photos_tikz/numerical_tests/Data/tripod_net/}
}{
}

\newcommand{\subfile}{0}
\centering
	\begin{tikzpicture}[baseline]
		\begin{axis}
		[	width=6cm,
			height=7cm,
			scale only axis,
			xmin = 0,
			xmax = 1,
			ymin = 0.9,
			ymax = 3.1,
			title = {1st edge},
			xlabel = {$x$},
			ylabel = {$\rho$},
			name=Axis1,
			legend cell align=left,
			legend style={at={(0,-0.25)},
						  anchor=north west,
						  },			  
							]

			\addplot+ [mark repeat = 50, mark phase = 47, thick,smooth] table[x = x,y = r_kinetic] {\pfad edge_1st\subfile.txt};
			\addlegendentry{Kinetic};

			\addplot+ [mark repeat = 49, mark phase = 46, thick,smooth] table[x = x,y = r_half] {\pfad edge_1st\subfile.txt};
			\addlegendentry{Half-moment};
			
			\addplot+ [mark repeat = 49, mark phase = 44, thick,smooth] table[x = x,y = r_p1] {\pfad edge_1st\subfile.txt};
			\addlegendentry{p1};
			\addplot+ [mark repeat = 49, mark phase = 47, thick,smooth] table[x = x,y = r_fag] {\pfad edge_1st\subfile.txt};
			\addlegendentry{p1\cite{2014arXiv1411.6109R}};
			
			\addplot+ [mark repeat = 49, mark phase = 47, thick,smooth] table[x = x,y = r_con] {\pfad edge_1st\subfile.txt};
			\addlegendentry{p1$(\rho_i=\rho_j)$};
			\addplot+ [mark repeat = 50, mark phase = 25, thick,smooth] table[x = x,y = r_keller] {\pfad edge_1st\subfile.txt};
			\addlegendentry{Keller-Segel};
			\addplot+ [mark repeat = 49, mark phase = 44, thick,smooth] table[x = x,y = r_in1] {\pfad edge_1st\subfile.txt};
			\addlegendentry{Initial condition};

		\end{axis}

		\begin{axis}[
			width=4cm,
			height=7cm,
			scale only axis,
			xmin = 0,
			xmax = 1,
			ymin = 0.9,
			ymax = 3.1,		
			xshift=0.8cm,
			at={(Axis1.east)},
			anchor=west,
			name=Axis2,					
			title = {3rd edge},
			xlabel = $x$,
			]
							]

			\addplot+ [mark repeat = 60, mark phase = 5, thick,smooth] table[x = x,y = r_kinetic] {\pfad edge_3rd\subfile.txt};

			\addplot+ [mark repeat = 60, mark phase = 5, thick,smooth] table[x = x,y = r_half] {\pfad edge_3rd\subfile.txt};
			\addplot+ [mark repeat = 60, mark phase = 2, thick,smooth] table[x = x,y = r_p1] {\pfad edge_3rd\subfile.txt};
			\addplot+ [mark repeat = 60, mark phase = 2, thick,smooth] table[x = x,y = r_fag] {\pfad edge_3rd\subfile.txt};
			
			\addplot+ [mark repeat = 60, mark phase = 4, thick,smooth] table[x = x,y = r_con] {\pfad edge_3rd\subfile.txt};
			\addplot+ [mark repeat = 60, mark phase = 20, thick,smooth] table[x = x,y = r_keller] {\pfad edge_3rd\subfile.txt};
			\addplot+ [mark repeat = 60, mark phase = 5, thick,smooth] table[x = x,y = r_in3] {\pfad edge_3rd\subfile.txt};
			
		\end{axis}
		
		\begin{axis}[
			width=4cm,
			height=4cm,
			scale only axis,
			xmin = 0,
			xmax = 1,
			yshift=-2cm,
			at={(Axis2.south)},
			anchor=north,
			name=Axis3,
			title = {2nd edge},
			xlabel = {$x$},
			ylabel = {$\rho$},
							]

			\addplot+ [mark repeat = 60, mark phase = 1, thick,smooth] table[x = x,y = r_kinetic] {\pfad edge_2nd\subfile.txt};

			\addplot+ [mark repeat = 60, mark phase = 1, thick,smooth] table[x = x,y = r_half] {\pfad edge_2nd\subfile.txt};
			\addplot+ [mark repeat = 60, mark phase = 1, thick,smooth] table[x = x,y = r_p1] {\pfad edge_2nd\subfile.txt};
			\addplot+ [mark repeat = 60, mark phase = 1, thick,smooth] table[x = x,y = r_fag] {\pfad edge_2nd\subfile.txt};
			
			\addplot+ [mark repeat = 60, mark phase = 1, thick,smooth] table[x = x,y = r_con] {\pfad edge_2nd\subfile.txt};
			\addplot+ [mark repeat = 60, mark phase = 1, thick,smooth] table[x = x,y = r_keller] {\pfad edge_2nd\subfile.txt};
			\addplot+ [mark repeat = 60, mark phase = 1, thick,smooth] table[x = x,y = r_in2] {\pfad edge_2nd\subfile.txt};
			
		\end{axis}
		\node[anchor=center] at (Axis1.center) {\large$\epsilon\approx 0$};

%
	\end{tikzpicture}%
\let\pfad\undefined
\caption{Numerical solutions on a tripod network at time $t=0.3$, $\Delta x=0.02$, $\epsilon = 10^-6$}
\label{fig:tripod_e0}
\end{figure}

In the figures \ref{fig:tripod_e1}-\ref{fig:tripod_e0} the densities at $t=0.3$ for different values of $\epsilon$ are shown.
For $\epsilon=1$ the solutions are similar to those on a single interval.
The half moment model is close to the kinetic one.
In case of the Cattaneo equations the solution using the coupling conditions \eqref{eq:cc_P1} is a good approximation of the kinetic model, whereas the continuity of $\rho$ is even worse than the Keller-Segel model.

As the value of $\epsilon$ decreases all solutions approach the state of the Keller-Segel equation.
The small deviations for $\epsilon=10^{-6}$ are caused by the different discretizations due to the relaxation scheme and disappear with refining the spatial resolution.

\subsection{Numerical solutions on a larger network}
\begin{figure}[ht!]
	\begin{subfigure}[b]{0.49\linewidth}
		\includegraphics[width=0.99\linewidth]{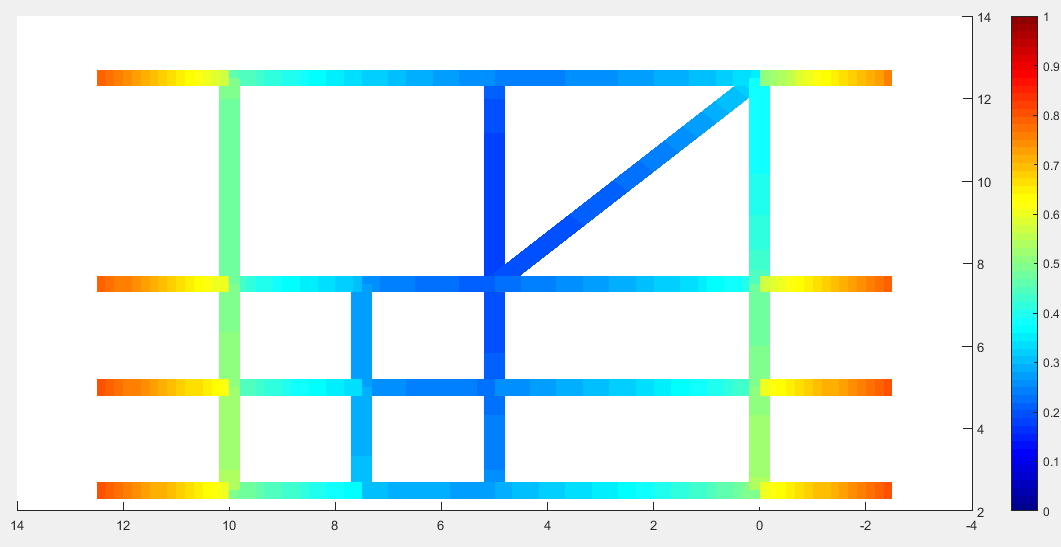}
		\caption{Kinetic model}
	\end{subfigure}
	\begin{subfigure}[b]{0.49\linewidth}
		\includegraphics[width=0.99\linewidth]{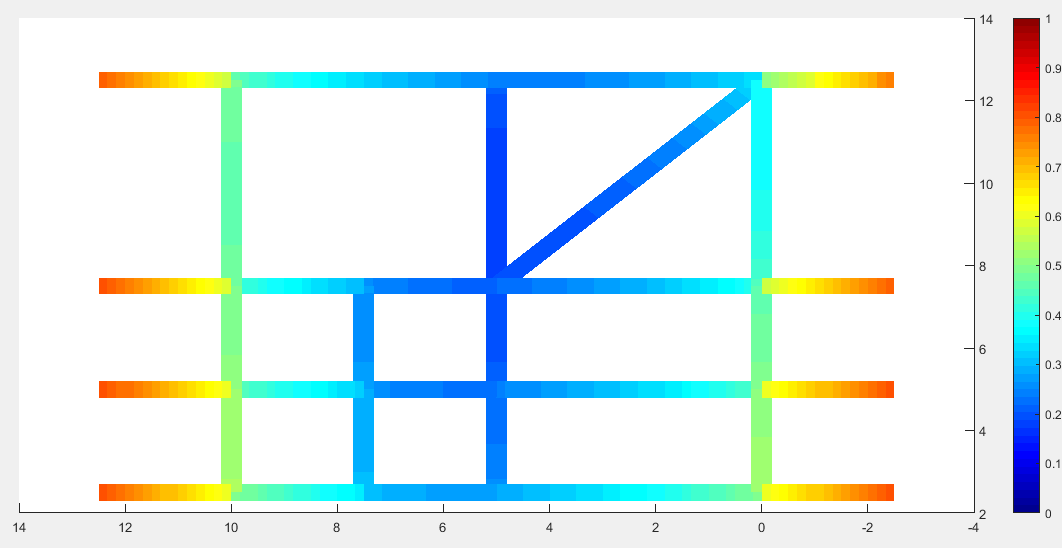}
		\caption{Half moment model}
	\end{subfigure}
	\\
	\begin{subfigure}[b]{0.49\linewidth}
		\includegraphics[width=0.99\linewidth]{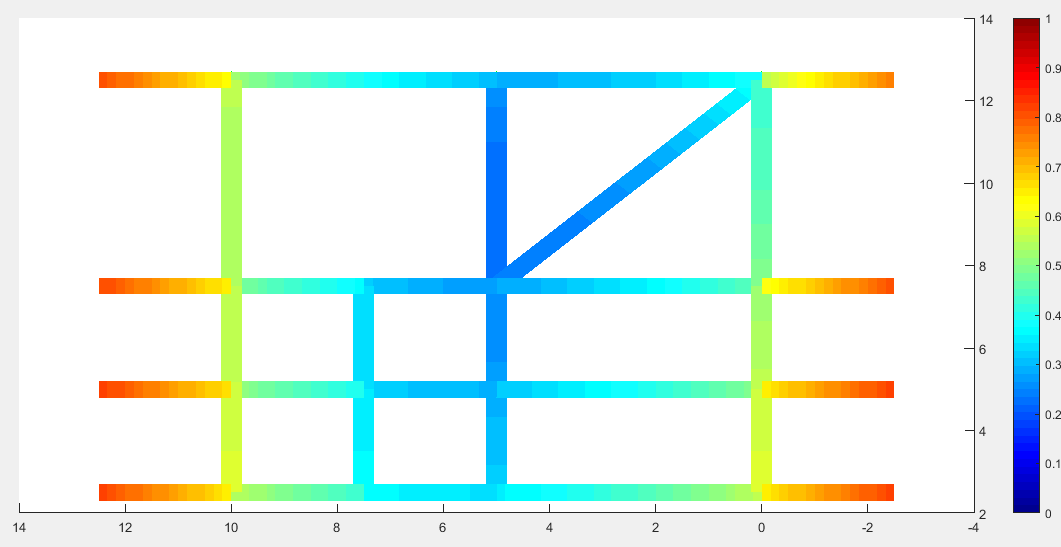}
		\caption{P1 model}
	\end{subfigure}
	\begin{subfigure}[b]{0.49\linewidth}
	\includegraphics[width=0.99\linewidth]{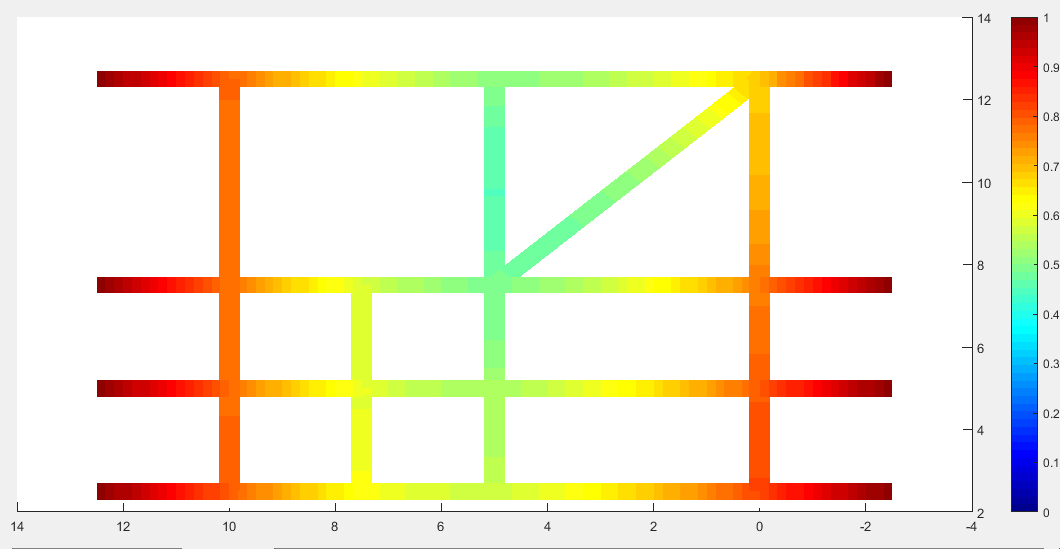}
		\caption{Keller-Segel model}
	\end{subfigure}
	
	\caption{Comparison of the numerical solutions of the four models on a larger network. at $t=5$.}
	\label{fig:large_net_t5}
\end{figure}

In this last test case we consider a larger network of $31$ edges and $23$ nodes as shown in Figure \ref{fig:large_net_t5}.
The length of the short edges is $0.5$, for the longer ones we have $1$ and $\sqrt{2}$ respectively.
Note that this network is not only build by three way junctions, but contains nodes connecting up to five edges.
At the open ends of the network Dirichlet boundary conditions are imposed.
For the kinetic model we prescribe $f(0,v,t)=\frac12$ for $v>0$. 
The boundary values for the other models can be derived thereof with $\epsilon = 0$. 
As initial conditions all values are set to zero except the density in the edges at the outer boundaries, which is set to $1$.
For the numerical computation we used 
$15$, $30$ and  $42$ cells for the spatial discretizaiton of the edges, respectively.
The velocity space in the kinetic model is resolved with $20$ points. 
In Figure \ref{fig:large_net_t5} the density at time $t=5$ for all four models is shown.
There is an inflow from both sides of the network.
As observed in the previous tests, the states of the Keller-Segel model propagate faster, such that the network is already filled by half.
The  solution of the kinetic model and the one of the half moment model almost coincide.
\begin{figure}[ht!]
	\begin{subfigure}[b]{0.49\linewidth}
		\includegraphics[width=0.99\linewidth]{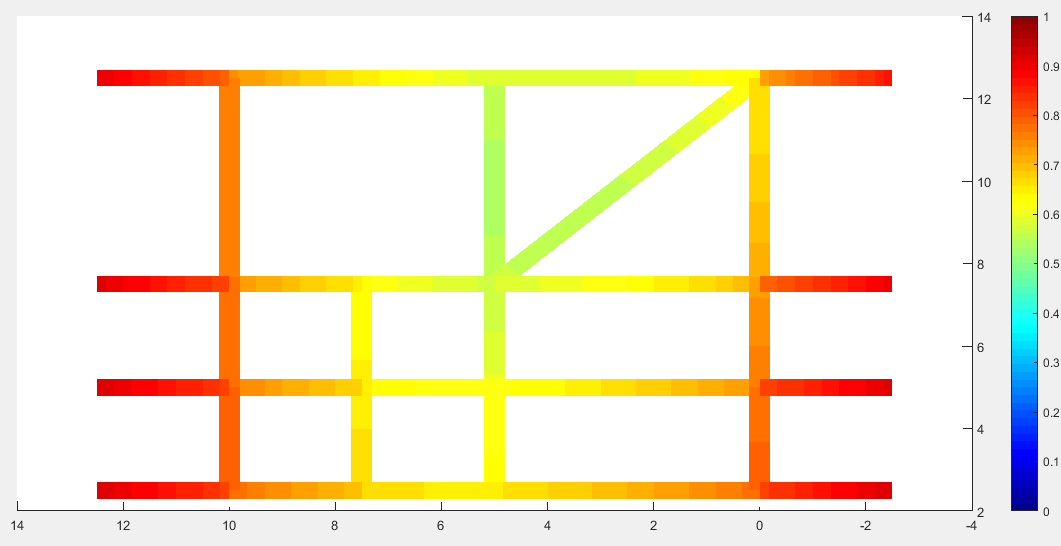}
		\caption{Kinetic model}
	\end{subfigure}
	\begin{subfigure}[b]{0.49\linewidth}
		\includegraphics[width=0.99\linewidth]{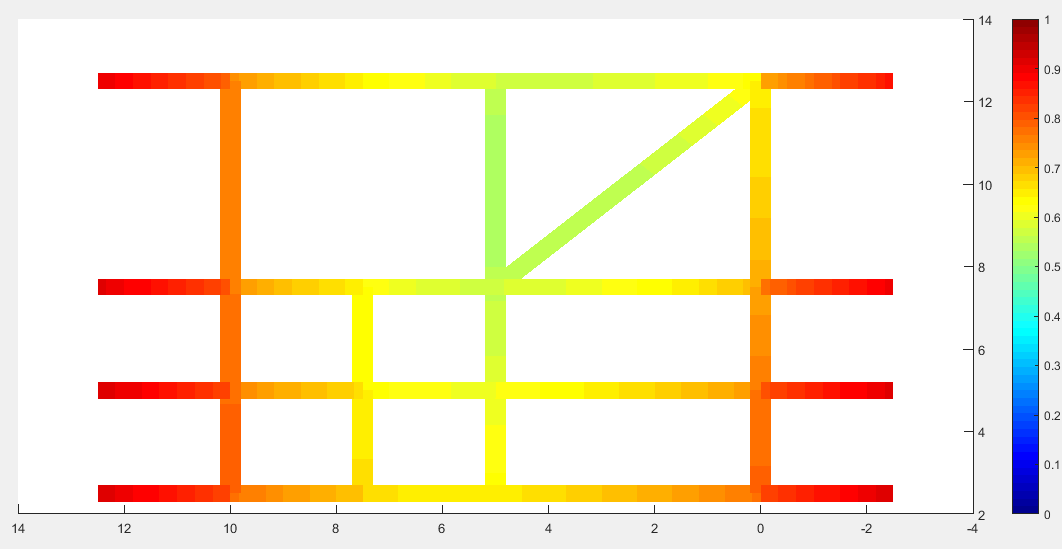}
		\caption{Half moment model}
	\end{subfigure}
	\\	
	\begin{subfigure}[b]{0.49\linewidth}
		\includegraphics[width=0.99\linewidth]{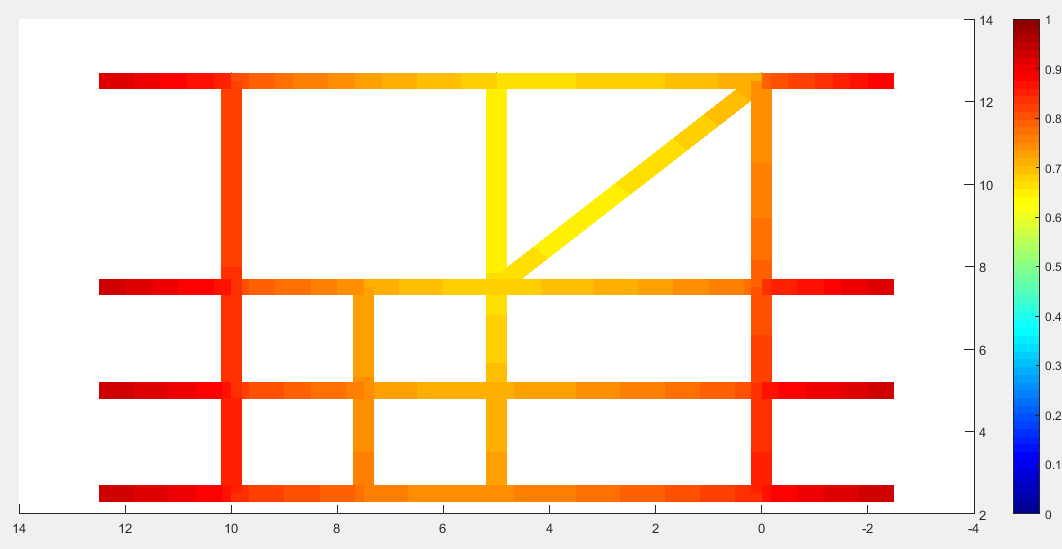}
		\caption{P1 model}
	\end{subfigure}
	\begin{subfigure}[b]{0.49\linewidth}
		\includegraphics[width=0.99\linewidth]{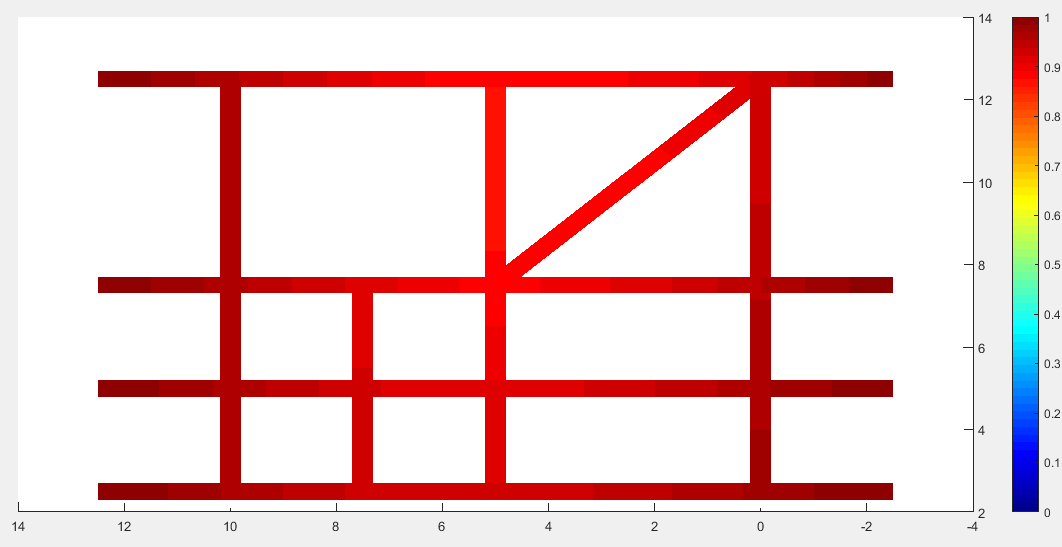}
		\caption{Keller-Segel model}
	\end{subfigure}
	\caption{Comparison of the numerical solutions on a larger network. at $t=15$.}
		\label{fig:large_net_t15}
\end{figure}
In Figure \ref{fig:large_net_t15} we can see the density at time $t=15$.
For the Keller Segel model the density has reached $1$ almost everywhere
and also the density in the P1 model is ahead of the kinetic and half moment model.
Note that in the lower part of the network the edges are filled more than in the upper part, as we have a denser grid and more inflow ends.
\begin{figure}[ht!]
\centering
\externaltikz{figure_mass_net}{
	\ifthenelse{\isundefined{\pfad}}{
    \def\pfad{./Photos_tikz/numerical_tests/Data/large_net/}
    }{
}

\begin{tikzpicture}
\begin{groupplot}[group style={group size=2 by 1, horizontal sep = 2cm,  vertical sep = 2cm},
width=10cm,
height=4cm,
scale only axis,
xmin = 0,
xmax = 30,
ymin = 0,
legend cell align=left,
legend style={at={(1,0.5)},
						  },			  
							]
\nextgroupplot[
title = {$\epsilon = 1$},
xlabel = {$t$},
ylabel = {mass of all network},
legend cell align=left,
]
\addplot+ [mark repeat = 10, mark phase = 5, smooth,line width=2pt] table[x = T,y = mass_net_kinetic] {\pfad kinetic/Mass_kinetic_t30_nx30.txt};
\addlegendentry{Kinetic};

\addplot+ [mark repeat = 10, mark phase = 10, smooth,line width=2pt] table[x = T,y = mass_net_half] {\pfad half_moment/Mass_half_t30_nx30.txt};
\addlegendentry{Half-moment};

\addplot+ [mark repeat = 10, mark phase = 15, smooth,line width=2pt] table[x = T,y = mass_net_full] {\pfad full_moment/Mass_full_t30_nx30.txt};
\addlegendentry{Full-moment};

\addplot+ [mark repeat = 10, mark phase = 20, smooth,line width=2pt] table[x = T,y = mass_net_keller] {\pfad keller_segel/Mass_keller_t30_nx30.txt};
\addlegendentry{Keller-Segel};

\end{groupplot}
\end{tikzpicture}
}
\caption{Total mass the large network.}
\label{fig:mass_net}
\end{figure}
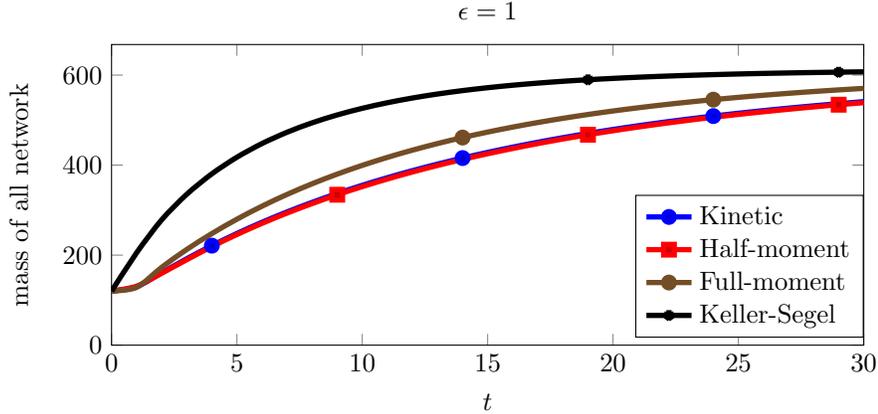
In Figure \ref{fig:mass_net} the evolution of the total density in the network up to $T=30$ is shown.
As before, we observe that the Keller-Segel model fills the network much faster, than the other three. 
The values of the half moment model and the kinetic one almost coincide. 
Thus the half moment model is a very good approximation to the kinetic one in this test case.

\section{Conclusions}
In this paper we developed a hierarchy of coupling conditions for chemotaxis models. 
A set of simple assumptions led to coupling conditions for the kinetic models, from which we derived corresponding equations for the averaged models. 
For more sophisticated dynamics on the network some conditions could be relaxed, e.g. nonlinear coupling relations, without affecting the presented approach.
This procedure also might be applied to other models, as the assumptions on the coupling are not specific to chemotaxis models.
In the numerical tests we investigated the dynamics for different coupling conditions and for varying values of $\epsilon$.
A simulation on a larger network showed the applicability to complicated structures and differences in behavior of macroscopic and kinetic  models.


\bibliographystyle{plain}
\bibliography{chemotaxis_coupling,lit.bib,lit2.bib}

\end{document}